\documentclass[a4paper]{amsart}

\usepackage[utf8]{inputenc}
\usepackage{amssymb,amsthm,amsmath}
\usepackage{hyperref}
\usepackage{enumitem}
\usepackage[nobysame]{amsrefs}
\usepackage{tikz}
\usepackage{xcolor}

\newtheorem{theorem}{Theorem}[section]
\newtheorem{proposition}[theorem]{Proposition}

\newtheorem{lemma}[theorem]{Lemma}
\newtheorem{corollary}[theorem]{Corollary}

\theoremstyle{definition}
\newtheorem{definition}[theorem]{Definition}
\newtheorem{example}[theorem]{Example}

\theoremstyle{remark}
\newtheorem{remark}[theorem]{Remark}

\newcommand{\dis}{\operatorname{dis}}
\newcommand{\Th}{\operatorname{Th}}
\newcommand{\join}{\triangledown}
\newcommand{\crossing}{\square}

\newcommand{\R}{\operatorname{R}}
\newcommand{\NR}{\operatorname{NR}}
\newcommand{\WR}{\operatorname{WR}}
\newcommand{\NWR}{\operatorname{NWR}}

\renewcommand{\emptyset}{\varnothing}

\title{Infinite atomized semilattices}
\author[]{Fernando Martin-Maroto}
\author[]{Antonio Ricciardo}
\thanks{Fernando Martin-Maroto and Antonio Ricciardo are joint first authors of this paper.}
\author[]{David Méndez}
\author[]{Gonzalo G. de Polavieja}
\address{Addresses of authors 1,3,4: Champalimaud Research, Champalimaud Foundation, Lisbon, Portugal}
\address{Addresses of authors 1,2,4: Algebraic AI, Madrid, Spain}


\begin{document}

\begin{abstract}
Atomized semilattices are subdirect product representations that facilitate the construction and manipulation of semilattices in practical contexts. This paper develops a careful axiomatization of atomized semilattices that extends to the infinite setting. We extend the Full Crossing operator, which yields the quotient via principal congruences, to infinite structures, and prove that every semilattice admits an atomization. The notions of atom redundancy are generalized to the infinite case, revealing a divergence between redundancy and weak redundancy that is absent in finite models. We also show that every semilattice has a freest completion, atomized by all the atoms compatible with the semilattice, and that every finitely generated semilattice admits an atomization consisting exclusively of non-redundant atoms.
\end{abstract}

\date{\today}
\maketitle

\tableofcontents

\section{Introduction}

In \cite{MarPol18}, the authors introduce Algebraic Machine Learning, a novel approach to learning from data by embedding data and rules and extracting rules from the data in an algebra, using the language of Universal Algebra, \cite{BurSan81}, and Model Theory, \cite{Mar02}.\footnote{A later account of Algebraic Machine Learning, which describes learning as computing subsets of a subdirect decomposition, appeared after the first version of this paper was written, \cite{MarAbdMenPol26}.} In particular, the algebras of choice are (join) semilattices. These structures are used due to them being simple enough that they are easy to work with, yet powerful enough that information of different kinds can be encoded in them using semantic embeddings, that is, embeddings that assume a decoder with no knowledge beyond the rules of the algebra, \cite{MarPol22}.

One key challenge on using these ideas lies on how to efficiently and effectively modify the algebras at work as to encode the logical rules and data that we are working with. For that purpose, the theory of atomized semilattices was developed, \cite{MarPol21}. Like in the case of Boolean algebras, \cite{Sto36}, and distributive lattices, \cite{Sto38}, we show that elements in a semilattice can be represented with sets that can be used to perform various computations. This is achieved using \emph{atomizations}.

An atomization on a semilattice consists of a collection of added elements, called atoms, that encode the order relation of the semilattice in a way that computations can be done exclusively via set-theoretic operations. As such, an atomization can also be regarded as a different representation of the semilattice structure. In particular, elements are characterized by the atoms lower than them, and two elements are comparable if and only if the sets of atoms lower than them are. Atoms are characterized by the set of constants greater than them, the constants being a collection of elements that generate the semilattice when taking their (formal) joins.

This theory is introduced in \cite{MarPol21} in the context of finite semilattices, and thus, it becomes natural to ask whether this can be extended to the infinite setting. Not only that, but it would also have computational implications in the framework of Algebraic Machine Learning, allowing for example for the representation of real values via infinite chains. The main objective of this paper is to settle that question. We do so following a structure similar to the one they used.

In Section \ref{section:preliminaries} we introduce the axiomatization of (possibly infinite) atomized semilattices (Definition \ref{definition:atomized_semilattice}) and show that it is still sound in the infinite setting. We assume that elements on a semilattice are generated as joins of a generating subset, which we call constants, and discuss how we keep track of them across algebras of terms (Definition \ref{definition:admissible_terms}), and how our axioms allow us to identify atoms across every model generated by the same constants. At the end of section \ref{section:preliminaries} we discuss the correspondence between atomizations of a semilattice and subdirect product decompositions into two-element semilattices whose factor maps preserve the joins named by the terms (Theorem \ref{theorem:subdirect}).

In Section \ref{section:full_crossing} we introduce the concept of relative freedom of models (Definition \ref{definition:freer_or_as_free}) together with the full crossing operation (Definitions \ref{definition:full_crossing} and \ref{definition:full_crossing_infinite}).\footnote{In later work we found that full crossing has a counterpart law in the lattice of congruences of a semilattice: the pairwise distributive law for maximal congruences, which is equivalent to crossing a single duple, \cite{MarRicPol26}*{Theorem 4.2 and Remark 4.5}.} We show that the full crossing of a duple in a model, which only requires performing set-theoretical operations on the atoms of the starting model, allows us to obtain the freest model that satisfies the crossed duple together with the positive theory of the starting model (Theorem \ref{theorem:full_crossing_freest} and Proposition \ref{proposition:full_crossing_freest}). As a consequence, the full crossing of a duple is equivalent to enforcing a principal congruence, and that of a finite set of duples to enforcing the congruence they generate, while for an infinite set of duples the enforced congruence may be larger (Theorem \ref{theorem:congruences}). Along the way, we also show that every semilattice admits an atomization and how they can be built (Theorem \ref{theorem:every_semilattice_atomizable}).

Finally, in Section \ref{section:atom_redundancy} we study atom redundancy. Unlike in the finite case, we will require two different notions of redundancy: redundancy and weak redundancy (Definition \ref{definition:redundancy}). In the case of complete semilattices with every non-empty set of constants as a term (which includes semilattices over a finite set of constants), and in the case of semilattices with no infinite strictly ascending chain of terms (Corollary \ref{corollary:ascending_chain_least_atomization}), they are equivalent concepts, but not in the general case. Nonetheless, we are able to show that in the finitely generated case, semilattices admit atomizations consisting exclusively of non-redundant atoms (Theorem \ref{theorem:finitely_generated_join_redundant} and Proposition \ref{proposition:completion_finitely_generated_NR} for an infinite set of constants, and \cite{MarPol21} for a finite one). This result can also be derived from the theorem of Birkhoff and Frink on algebraic lattices, \cite{BirFri48} (Remark \ref{remark:Birkhoff_Frink}). We also show that the atoms common to all atomizations are the non-weakly-redundant ones (Corollary \ref{corollary:every_atomization_contains_non-weekly_redundant}), and that every semilattice has a freest completion, atomized by all the atoms compatible with it (Theorem \ref{theorem:freest_completion}).

\section{Preliminaries}\label{section:preliminaries}

This section follows closely and is a generalization of \cite{MarPol21}*{Section 2}. The objective here is to show that semilattice atomizations can be extended to the infinite setting.

Remember that a semilattice is an algebra with a single binary operation $\odot$ that is commutative, associative and idempotent. Equivalently, a (join) semilattice is a partially ordered set that has a least upper bound for every non-empty finite subset of elements. The translation between these two formulations is that $a\le b$ in a semilattice $M$ if and only if $M\models (b = a\odot b)$.

Since $\odot$ is binary, the only joins that are part of the semilattice structure are the finite ones. For an infinite family $\{a_i\}_{i\in I}$ of elements of $M$, the join $\odot_{i\in I} a_i$ means its least upper bound in the order of $M$, which may or may not exist. An infinite join is therefore not an operation of the semilattice signature: it is an order-theoretic supremum, present or absent depending on $M$. Accordingly, a semilattice homomorphism, being a map that preserves $\odot$, preserves finite joins but is not required to preserve infinite ones, even when they exist in both its domain and its codomain. For instance, in the chain $0<1<2<\cdots<\infty$ we have $\odot_{n\in\mathbb{N}}\, n = \infty$, yet the map onto the two-element semilattice $\{0<1\}$ sending every $n\in\mathbb{N}$ to $0$ and $\infty$ to $1$ is a semilattice homomorphism.

Given our interest in computational applications where we keep modifying the underlying algebra to model provided rules and data, we want to perform subsequent quotients of a semilattice while being able to refer to their elements by the same names. We can assume without loss of generality that there is a set of \emph{constants} $C$ and a map $C \to M$ such that every element of $M$ can be obtained as a (finite or infinite) join of elements in the image of $C$. Indeed, it is always possible to consider $C = M$ and the map to be the identity, although this is very rarely necessary. We say that $M$ is a semilattice over (the set of constants) $C$. 

We will mostly be interested in two different types of infinite semilattices.

\begin{definition}
    Let $M$ be a semilattice defined over a set of constants $C$. We say that $M$ is \emph{finitely generated} over $C$ if every element of $M$ can be written as a join of a finite number of elements of $C$.
\end{definition}

\begin{definition}
    Let $M$ be a semilattice defined over a set of constants $C$. We say that $M$ is \emph{complete} if the join of any non-empty subset of elements of $C$ is an element in $M$.
\end{definition}

We now introduce the theoretical background allowing us to refer to the elements of $M$ as the join of the involved constants. 

\begin{definition}
    Let $C$ be a (finite or infinite) set of constants. The freest (finitely generated) semilattice over $C$, denoted $F_C(\emptyset)$, is the semilattice that has one different element per each element of $C$, together with every finite idempotent summation of elements of $C$. Therefore, it is the set of all finite non-empty subsets of $C$.
\end{definition}

If $M$ is a semilattice for which there is a map $C\to M$ as introduced above, there is a homomorphism
\begin{equation*}
    \nu_M \colon F_C(\emptyset) \longrightarrow M
\end{equation*}
extending it. We then say that $M$ is a semilattice over $C$.

The homomorphism $\nu_M$ is surjective if $M$ is finitely generated and, generally, it allows us to consider finitely generated terms over $C$ in any semilattice over $C$, by means of taking the homomorphic image of said term in $M$. It may be the case that in a model $M$, $M\models \big(\nu_M(t) = \nu_M(s)\big)$ for different terms $t$ and $s$, in which case we will denote $M\models (t=s)$.

Notice that every finite semilattice is both finitely generated and complete. Also notice that a semilattice $M$ being finitely generated over $C$ depends on the choice of $C$. In particular, if we choose $C = M$, it will be finitely generated over $C$. However, a semilattice being complete is independent of the choice of $C$. Indeed, we are requiring that every join of elements of $C$ can be formed, which must also include every possible join of elements of $M$.


\begin{definition}
    Let $C$ be a (finite or infinite) set of constants. The freest complete semilattice over $C$, denoted $\tilde{F}_C(\emptyset)$, is the semilattice that has one different element per each element of $C$, together with every (finite and infinite) idempotent summation of elements of $C$. Therefore, it is the set of all non-empty subsets of $C$.
\end{definition}

If $M$ is a complete semilattice, we are also guaranteed the existence of a homomorphism
\begin{equation*}
    \tilde{\nu}_M \colon \tilde{F}_C(\emptyset) \longrightarrow M.
\end{equation*}
As such, we can use the freest complete semilattice to keep track of the elements of $M$. However, such homomorphism may not exist if $M$ is not complete. We thus need a more general way of introducing this homomorphism, based on the joins of constants that exist on a particular model.

\begin{definition}
    Let $C$ be a set of constants and let $T$ be a set of non-empty subsets of $C$ such that $\{c\}\in T$, $\forall c\in C$. Assume that $T$ is closed under finite unions, thus it is a semilattice. We can view it as a semilattice over $C$ where the join is the set union. We denote it by $F_{C;T}(\emptyset)$, and we say $F_{C;T}(\emptyset)$ is the \emph{freest semilattice over $C$ with terms $T$}. We say that $T$ is a \emph{set of terms}, and an element $t\in T$ is implicitly seen as the element $t\in F_{C;T}(\emptyset)$. We then use the notation $C(t)$ to refer to the set of constants mentioned in $t$.
\end{definition}

And thus, we can introduce the following.

\begin{definition}\label{definition:admissible_terms}
    Let $M$ be a semilattice over a set of constants $C$. Let $T$ be a set of non-empty subsets of $C$ closed under finite unions such that $\{c\}\in T$, for all $c\in C$. We say that $T$ is a \emph{set of terms} for $M$ if every element of $M$ can be expressed as the join of constants in $t$ for some $t\in T$; and for every $t\in T$, the join of the constants in $t$ is well-defined in $M$. The elements of $T$ are then referred to as \emph{admissible terms} or just \emph{terms} of $M$.
\end{definition}

Let $M$ be a semilattice over a constant set $C$ with terms $T$. Since $T$ is closed under finite unions, $T$ is itself a semilattice over $C$ where the join is the set union. We then have that $F_C(\emptyset)\subseteq F_{C;T}(\emptyset)\subseteq \tilde{F}_C(\emptyset)$. And, since the joins of elements of $T$ exist in $M$, there is always a homomorphism
\[\nu_{M;T}\colon F_{C;T}(\emptyset)\longrightarrow M,\]
which is also surjective. We will often denote $\nu_{M;T}=\nu_M$, since the value $\nu_{M;T}(t)$, the join in $M$ of the constants of $t$, does not depend on $T$.

We remark that more than one valid set of terms $T$ may exist for a semilattice $M$, and the choice of $T$ will affect some constructions carried out throughout this paper. The set $T$ should be thought of a selection of which terms we allow to be represented, although they will always be enough so that every element of $M$ can be attained. 

In particular, if $M$ is finitely generated, we can always choose $T$ to be the set of finite subsets of $C$, in which case $F_{C;T}(\emptyset) = F_C(\emptyset)$. Generally, whenever we are talking about finitely generated semilattices, we are implicitly making this selection for $T$. On the other hand, if $M$ is complete, we can choose $T$ to be the power set of $C$ (excluding the empty set), thus $F_{C;T}(\emptyset) = \tilde{F}_C(\emptyset)$. Notice further that, no matter the choice of $T$, there is an injective homomorphism $F_{C;T}(\emptyset)\to \tilde{F}_C(\emptyset)$, which lets us consider the terms of $M$ as terms in the freest complete semilattice even when $M$ is not complete.

\begin{definition}\label{definition:duples}
    A \emph{duple} over a constant set $C$ is an ordered pair of elements $r=(r_L, r_R)$ of $\tilde{F}_C(\emptyset)$. If $M$ is a semilattice over $C$ with terms $T$, a \emph{duple of $M$} is an ordered pair $r=(r_L, r_R)$ of elements of $F_{C;T}(\emptyset)$. We consider positive and negative duples, $r^+$ and $r^-$, where we say $M\models r^+$ if $M\models \big(\nu_{M}(r_L)\le \nu_{M}(r_R)\big)$; and $M\models r^-$ otherwise. If $r_L, r_R\in F_C(\emptyset)$ we say that $r$ is a \emph{finitely generated} duple. Note that if $M$ is any semilattice over $C$ with terms $T$ and $r$ is a finitely generated duple, it is always the case that $r_L, r_R\in F_{C;T}(\emptyset)$, so $r$ is a duple of $M$.
\end{definition}

Maps between semilattices over $C$ should respect not only the constants but also the terms.

\begin{definition}\label{definition:interpretation_preserving}
    Let $M$ and $N$ be semilattices over a set of constants $C$, and let $T$ and $T'$ be sets of terms for $M$ and $N$ respectively, with $T\subseteq T'$. A semilattice homomorphism $f\colon M\to N$ is \emph{interpretation-preserving} if $f(\nu_M(t)) = \nu_N(t)$ for every $t\in T$.
\end{definition}

Taking $t=\{c\}$ gives $f(c)=c$ for every $c\in C$. Since $\nu_M(t)$ is the join of the constants of $t$, the condition reads $f\big(\odot_{c\in C(t)}\, c\big) = \odot_{c\in C(t)} f(c)$ for every $t\in T$. Thus an interpretation-preserving homomorphism preserves, besides the finite joins, the infinite joins named by $T$, extending the notion of homomorphism to the infinite joins we have chosen to represent. As $f$ is order-preserving, $f(\nu_M(t))\ge\nu_N(t)$ always holds, and only the reverse inequality may fail. If every term of $T$ is finite, as for a finitely generated $M$, interpretation-preserving simply means constant-preserving, $f(c)=c$ for $c\in C$. For infinite terms the condition is stronger. In the chain $0<1<2<\cdots<\infty$ above, seen as a semilattice over $C=\{0,1,2,\dots,\infty\}$ with all non-empty subsets of $C$ as terms, consider $\{0<1\}$ as a semilattice over $C$ through the map sending every $n\in\mathbb{N}$ to $0$ and $\infty$ to $1$. This map preserves the constants, but not the interpretation of the term $t$ with $C(t)=\mathbb{N}$, since $\nu(t)=\infty$ is sent to $1$ while $\odot_{n\in\mathbb{N}}\, 0 = 0$.

The following concept will also be useful.

\begin{definition}\label{definition:completion}
    Let $S$ be a semilattice over a constant set $C$ with terms $T$. We say that a complete semilattice $M$ over $C$ is a \emph{completion of $S$} if there is an injective, interpretation-preserving semilattice homomorphism $i\colon S\to M$, Definition~\ref{definition:interpretation_preserving}, where $M$ has as terms all non-empty subsets of $C$. 
\end{definition}

With this, we are now ready to define atomized semilattices.

\begin{definition}\label{definition:atomized_semilattice}
    An \emph{atomized semilattice} $M$ over a set of constants $C$ with terms $T$ is a semilattice $(M,\odot)$ over $C$ with terms $T$ (whose elements are called \emph{regular elements} and are denoted by Latin letters) together with a set of special elements, called \emph{atoms} (denoted with Greek letters), such that the partial order relation in $M$ can be extended to include the atoms and such that the following axioms are satisfied.
    \begin{enumerate}[label=(AS\arabic*), ref=AS\arabic*]
        \item \label{axiom:AS1} $\forall\phi\exists c\colon (c\in C)\wedge (\phi < c)$;
        \item \label{axiom:AS2} $\forall\phi\forall a$, $(a\not\le\phi)$;
        \item \label{axiom:AS3} $\forall a\forall b \big(a\le b \Leftrightarrow \neg\exists\phi\colon ((\phi < a) \wedge (\phi\nless b))\big)$;
        \item \label{axiom:AS4} $\forall\phi\forall t\in T \big((\phi < \nu_M(t)) \Leftrightarrow (\exists c\in C(t)\colon\phi < \nu_M(c))\big)$;
        \item \label{axiom:AS5} $\forall\phi\forall\psi\Big(\forall c\in C \big((\phi < c) \Leftrightarrow (\psi<c)\big)\Rightarrow (\phi=\psi)\Big)$;
        \item \label{axiom:AS6} $\forall a\exists \phi\colon (\phi < a)$.
    \end{enumerate}
    We also introduce the following associated notation.
    \begin{itemize}
        \item If $b\in M$ is a regular element, $L_M^a(b)$ is the \emph{lower atomic segment} of $b$, that is, the set of atoms $\phi$ satisfying $M\models (\phi < b)$;
        \item the \emph{upper constant segment} $U^c(\phi)$ of an atom $\phi$ is the set of constants $c\in C$ such that $M\models(\phi < c)$.
    \end{itemize}
    Notice that in this definition we are identifying $C$ with $\{\nu_M(c)\mid c\in C\}$. The reason why we use $M$ as a subindex in the lower atomic segment but not in the upper constant segment is that atoms are defined universally, independent of the model by their upper constant segment. However, the atoms in the lower segment of a term depend upon the model. Notice further that $U^c(\phi)$ may be an infinite set, and a semilattice may have an infinite amount of atoms.  We introduce the notation $\phi:< U$ to refer to an atom whose upper constant segment is $U$, $U^c(\phi) = U$, which by \eqref{axiom:AS5} is enough to fully characterize the atom. Finally, notice that it may be the case that multiple terms give rise to equal regular elements.
\end{definition}

\begin{remark}
    Using this notation, the order relation in $F_{C;T}(\emptyset)$ is characterized by the fact that given $t, s\in T$, $F_{C;T}(\emptyset)\models t\le s$ if and only if $C(t)\subseteq C(s)$. Indeed, $F_{C;T}(\emptyset)\models t\le s$ if and only if $s = t\odot s$, but in $T$, the semilattice operation is the union. Thus, $s = t\odot s$ is equivalent to saying that $C(s) = C(s)\cup C(t)$, which is itself equivalent to $C(t)\subseteq C(s)$.
\end{remark}

We will quite often abuse notation and regard a term $t\in T$ as an element of $M$, implicitly referring to $\nu_{M}(t)$.

\begin{remark}
    In \cite{MarPol21}, \eqref{axiom:AS4} is formulated as
    \begin{equation}\label{axiom:AS4f}
        \tag{AS4f} \forall\phi\forall a\forall b \big(\phi<a\odot b\Leftrightarrow (\phi < a)\vee (\phi < b)\big).
    \end{equation}
    Both formulations are equivalent when every term of $T$ is finite, in particular in the finite case, since \eqref{axiom:AS4} then follows from applying \eqref{axiom:AS4f} finitely many times. With infinite terms, \eqref{axiom:AS4f} alone would not determine whether an atom is lower than an infinite join of constants none of which is in its upper constant segment, which is why \eqref{axiom:AS4} is stated for whole terms.

    Notice that \eqref{axiom:AS4} quantifies over the terms of $T$, not over every family whose join happens to exist in $M$. As recalled above, such a join is a supremum, not an operation, and atoms need not respect it. For instance, in the ascending chain $c_0<c_1<\cdots<c_\infty$ considered below, with $T$ the finite subsets of $C$, the join of the constants $c_i$ is $c_\infty$, yet the atom $\phi_\infty:<\{c_\infty\}$ is lower than $c_\infty$ and lower than no $c_i$. Consequently, which sets of atoms atomize $M$ depends on the choice of $T$.
\end{remark}

Consider, for instance, a semilattice $M$ over $C=\{c_0, c_1, \dots, c_\infty\}$ whose constants form the ascending chain $c_0<c_1<c_2<\cdots<c_\infty$, and the infinite term $t$ with $C(t)=\{c_i\mid i\in\mathbb{N}\}$. If $T$ is the set of finite subsets of $C$, so that $t\notin T$, the elements of $M$ are the values of finite terms, that is, $c_0, c_1, \dots$ and $c_\infty$. The only upper bound of the $c_i$ in $M$ is $c_\infty$, so their supremum in $M$ is $\odot_{i\in\mathbb{N}} c_i = c_\infty$. If $t\in T$, then $\nu_M(t)$ is itself an upper bound of the $c_i$, and $\odot_{i\in\mathbb{N}} c_i = \nu_M(t)\le c_\infty$, where the inequality is an equality if $M$ consists only of the constants, and is strict if $\nu_M(t)$ is a further element of $M$, above every $c_i$ and below $c_\infty$.

Just like in the finite case, the order relation induced on terms by the atomization coincides with the order relation pre-existing in the semilattice. Indeed, consider the equation
\begin{equation}\label{eq:AS3b}
    \tag{AS3b} \forall a\forall b (a\le b\Leftrightarrow a\odot b = b).
\end{equation}
Here, as in \cite{MarPol21}, $\le$ may be regarded as a partial order on regular elements and atoms that is not a priori tied to $\odot$, so that \eqref{eq:AS3b} states that its restriction to regular elements is the order of the semilattice, as assumed in Definition \ref{definition:atomized_semilattice}. It turns out that \eqref{eq:AS3b} is weaker but almost equivalent to \eqref{axiom:AS3}, in the way we now show.

\begin{proposition}
    Assume \eqref{axiom:AS4}. Then,
    \begin{enumerate}
        \item \eqref{axiom:AS3} implies \eqref{eq:AS3b}.
        \item If $\forall a\forall b\big(\forall\phi ((\phi < a) \Leftrightarrow (\phi < b))\Rightarrow (a = b)\big)$, then \eqref{eq:AS3b} implies \eqref{axiom:AS3}.
        \item \eqref{axiom:AS3} implies $\forall a\forall b\big(\forall\phi ((\phi < a) \Leftrightarrow (\phi < b))\Rightarrow (a = b)\big)$.
    \end{enumerate}
\end{proposition}

\begin{proof}
    We first show that \eqref{axiom:AS4} implies \eqref{axiom:AS4f}. Let $a,b\in M$. Since $T$ is a set of terms for $M$, there are $t,s\in T$ with $a=\nu_M(t)$ and $b=\nu_M(s)$. As $T$ is closed under finite unions, $t\cup s\in T$ and $\nu_M(t\cup s)=a\odot b$. By \eqref{axiom:AS4}, $\phi<a\odot b$ if and only if $\phi<\nu_M(c)$ for some $c\in C(t)\cup C(s)$, that is, if and only if $\phi<a$ or $\phi<b$. With \eqref{axiom:AS4f} at hand, the proof of \cite{MarPol21}*{Theorem 1} applies verbatim, as it does not make use of the finiteness of the semilattice.
\end{proof}

Notice that, unlike in the finite case, the terms admissible in $M$ are not fully determined from $C$, as there may be multiple valid choices for the set of terms $T$. Which sets of atoms atomize $M$ depends on this choice, since \eqref{axiom:AS4} is required for every term of $T$. If $T\subseteq T'$ are two sets of terms for $M$, every atomization of $M$ with respect to $T'$ is also an atomization of $M$ with respect to $T$, since \eqref{axiom:AS1}--\eqref{axiom:AS3}, \eqref{axiom:AS5} and \eqref{axiom:AS6} do not mention $T$, and \eqref{axiom:AS4} for $T'$ includes \eqref{axiom:AS4} for $T$. The converse may fail. For instance, in the ascending chain above, with $M$ consisting only of the constants, every subset of $C$ has a join in $M$, so both the set of finite subsets of $C$ and the set of all non-empty subsets of $C$ are sets of terms for $M$. The atoms $\phi_n:<\{c_n, c_{n+1},\dots,c_\infty\}$, $n\in\mathbb{N}$, together with $\phi_\infty:<\{c_\infty\}$, atomize $M$ with respect to the former but not the latter: for the term $t$ with $C(t)=\{c_i\mid i\in\mathbb{N}\}$ we have $\phi_\infty<c_\infty=\nu_M(t)$, while $\phi_\infty$ is lower than no constant in $C(t)$. We now introduce some results on compatibility between the notations we just introduced.

\begin{proposition}\label{proposition:atomizations_and_order}
    Let $C$ be a set of constants and let $M$ be an atomized semilattice over $C$ with terms $T$. Let $t\in F_{C;T}(\emptyset)$ be a term and consider $\nu_M(t)$. Let $\phi$ be an atom, $c\in C$ a constant and $a, b$ regular elements.
    \begin{enumerate}
        \item $\forall t\forall c \big(c\in C(t)\Rightarrow \nu_M(c)\le \nu_M(t)\big)$;
        \item $\phi < \nu_M(t)\Leftrightarrow\exists c\colon \big((c\in C(t))\wedge (\phi<\nu_M(c))\big)$;
        \item  $\forall a \forall t \big((a =  \nu_M(t))  \Rightarrow  L_M^a(a) = \{\phi \text{ an atom of } M\mid C(t)\cap U^c(\phi)\ne\emptyset\}\big)$;
        \item if $t_i\in T$, $i\in I$, are terms such that $\cup_{i\in I} t_i\in T$, then $\odot_{i\in I} \nu_M(t_i) = \nu_M(\cup_{i\in I} t_i)$ and $L_M^a(\odot_{i\in I} \nu_M(t_i)) = \cup_{i\in I} L_M^a(\nu_M(t_i))$;
        \item $a\le b\Leftrightarrow L_M^a(a)\subseteq L_M^a(b)$.
    \end{enumerate}
\end{proposition}

\begin{proof}
    We prove each of the points independently.
    \begin{enumerate}
        \item Since $c\in C(t)$, $t = t\odot c$. The result then follows from the fact that $\nu_M$ is a homomorphism.
        \item Assume first that $\exists c\colon \big((c\in C(t))\wedge (\phi<\nu_M(c))\big)$. Then from (1), $\phi < \nu_M(c)\le \nu_M(t)$, and by the transitive property of the order relation, $\phi<\nu_M(t)$. On the other hand, if $\phi<\nu_M(t)$, by applying \eqref{axiom:AS4} we obtain $\phi<\nu_M(c)$ for some $c\in C(t)$.

        \item Assume $a = \nu_M(t)$. Then, $\phi < a$ if and only if $\phi < \nu_M(t)$. Assume first that $\phi$ is an atom of $M$ such that $C(t)\cap U^c(\phi)\ne\emptyset$. Let $c\in C(t)\cap U^c(\phi)$. By (1), $\nu_M(c)\le \nu_M(t) = a$. Furthermore, $\phi < \nu_M(c)$ since $c\in U^c(\phi)$. Thus, by (2), $\phi < \nu_M(t) = a$, that is, $\phi\in L_M^a(a)$. Reciprocally, if $\phi\in L_M^a(a)$, $\phi <\nu_M(t) = a$, and by (2) there is $c\in C(t)$ with $\phi <\nu_M(c)$. Therefore, $c\in C(t)\cap U^c(\phi)\ne\emptyset$.

        \item Write $t = \cup_{i\in I} t_i\in T$. Since $\nu_M(t)$ is the join in $M$ of the constants in $C(t)=\cup_{i\in I} C(t_i)$, it is an upper bound of every $\nu_M(t_i)$. Any upper bound of all the $\nu_M(t_i)$ is an upper bound of all the constants in $C(t)$, hence greater than or equal to $\nu_M(t)$. So $\odot_{i\in I}\nu_M(t_i) = \nu_M(t)$. Now by (3), $\phi < \nu_M(t)$ if and only if $C(t)\cap U^c(\phi)\ne\emptyset$, if and only if $C(t_i)\cap U^c(\phi)\ne\emptyset$ for some $i\in I$, if and only if $\phi < \nu_M(t_i)$ for some $i\in I$.

        \item By \eqref{axiom:AS3}, $a\le b$ if and only if there is no atom $\phi$ with $\phi<a$ and $\phi\nless b$, that is, if and only if $L_M^a(a)\subseteq L_M^a(b)$. \qedhere
    \end{enumerate}
\end{proof}

Properties (3)--(5) exhibit that knowing the admissible terms in an atomized semilattice and the constants in the upper constant segment of each of its atoms is enough to have the semilattice completely characterized.

\begin{definition}\label{definition:atom_from_ucs}
    Let $U$ be a non-empty subset of a set of constants $C$. The \emph{atom over $C$ with upper constant segment $U$}, denoted $\phi:<U$, is the atom $\phi$ with $U^c(\phi) = U$. For a term $t$ over $C$ we write $\phi<t$ if $U\cap C(t)\ne\emptyset$, and $\phi\nless t$ otherwise. We say that $\phi$ \emph{discriminates} a pair of terms $(t, s)$ if $\phi<t$ and $\phi\nless s$.
\end{definition}

By \eqref{axiom:AS5}, an atomized semilattice over $C$ has at most one atom with upper constant segment $U$, and if $\phi$ belongs to an atomization of a semilattice $M$ with terms $T$, then $\phi<\nu_M(t)$ if and only if $\phi<t$, for every $t\in T$, by Proposition \ref{proposition:atomizations_and_order} (2). Thus, an atom can be described independently of any atomized semilattice.

\begin{definition}\label{definition:compatible}
    Let $M$ be a semilattice over $C$ with terms $T$. We say that an atom $\phi$ over $C$ is \emph{compatible with $M$} or just that $\phi$ is in $M$, written $\phi\in M$, if every duple of $M$ discriminated by $\phi$ is negative in $M$. We denote the set of all atoms compatible with $M$ by $\Omega(M)$.
\end{definition}

The definition does not require $M$ to be atomized. If $M$ is atomized by a set of atoms $A$, then $A\subseteq\Omega(M)$ by \eqref{axiom:AS3}. Furthermore, $A\cup\{\phi\}$ atomizes $M$ if and only if $\phi\in\Omega(M)$. Indeed, \eqref{axiom:AS1}, \eqref{axiom:AS2}, \eqref{axiom:AS4} and \eqref{axiom:AS5} hold for any set of atoms of $\Omega(M)$ with the reading of Definition \ref{definition:atom_from_ucs}, which is well defined for such atoms: if $\nu_M(t) = \nu_M(t')$, the duples $(t,t')$ and $(t',t)$ are positive, so a compatible atom discriminates neither and is lower than $t$ if and only if it is lower than $t'$. \eqref{axiom:AS6} holds for $A\cup\{\phi\}$ because it holds for $A$, and adding $\phi$ affects \eqref{axiom:AS3} only through the duples discriminated by $\phi$. Hence Definition \ref{definition:compatible} agrees with the definitions of \cite{MarPol21}*{Definitions 2.3 and 2.4}: $\phi\in M$ if and only if $A\cup\{\phi\}$ atomizes the same semilattice as $A$, if and only if some atomization of $M$ contains $\phi$. By the same argument, if $A$ atomizes $M$ and $A\subseteq B\subseteq\Omega(M)$, then $B$ atomizes $M$: \eqref{axiom:AS6} holds for $B$ because it holds for $A$, and \eqref{axiom:AS3} holds because the atoms of $B$ discriminate no positive duple, while those of $A$ already discriminate every negative one. In particular, $\Omega(M)$ atomizes $M$, since every semilattice is atomizable, Theorem \ref{theorem:every_semilattice_atomizable}. Another possible rewrite of \eqref{axiom:AS3} is then
\begin{equation*}
    (t\le_M s) \Longleftrightarrow \neg\exists\phi\colon\big((\phi\in M)\wedge (\phi < t)\wedge (\phi\nless s)\big).
\end{equation*}

\begin{proposition}\label{proposition:all_atoms_fit_freest}
    Let $C$ be a set of constants and let $T$ be a set of terms over $C$. If $A$ is a set of atoms that atomizes $F_{C;T}(\emptyset)$ and $\{\psi_i\mid i\in I\}$ is any set of atoms, then $A\cup \{\psi_i \mid i\in I\}$ is an atomization for $F_{C;T}(\emptyset)$.
\end{proposition}

\begin{proof}
    Every semilattice satisfies that if $C(t)\subseteq C(s)$ for two terms $t$ and $s$, then $t\le s$. The freest model is the only one for which the converse is also true. Since $A$ atomizes $F_{C;T}(\emptyset)$, axioms \eqref{axiom:AS1}, \eqref{axiom:AS5} and \eqref{axiom:AS6} are satisfied by $A$ and $A\cup \{\psi_i \mid i\in I\}$, where we are assuming $\psi\in A$ if an atom with upper constant segment $U^c(\psi)$ is already in $A$, following \eqref{axiom:AS5}. Furthermore, \eqref{axiom:AS2} and \eqref{axiom:AS4} hold by construction. Thus, we only need to check that $A\cup \{\psi_i \mid i\in I\}$ still satisfies \eqref{axiom:AS3}.

    Now, if $t\le s$, that is, if $C(t)\subseteq C(s)$, then $(\psi_i < t)\Rightarrow (\psi_i < s)$, for all $i\in I$, by Definition \ref{definition:atom_from_ucs}, since $U^c(\psi_i)\cap C(t)\subseteq U^c(\psi_i)\cap C(s)$. If, on the other hand, $C(t)\not\subseteq C(s)$, if $A$ satisfies \eqref{axiom:AS3}, $A\cup \{\psi_i \mid i\in I\}$  must also do so.
\end{proof}

This shows that if an atomization of the freest models exists, it will remain an atomization for the freest model after we add any number of atoms. In particular, the set of all atoms defined over $C$ would atomize $F_{C;T}(\emptyset)$. We now show that atomizations for these freest models do exist and can be as small as $|C|$.

\begin{proposition}\label{proposition:freest_admits_atomization}
    Let $C$ be a set of constants and let $T$ be a set of terms over $C$. Let $A = \{\psi_c :< \{c\}\mid c\in C\}$. Then $A$ atomizes $F_{C;T}(\emptyset)$.
\end{proposition}

\begin{proof}
    Remember that the freest model $F_{C;T}(\emptyset)$ is characterized by the property that $F_{C;T}(\emptyset) \models s\le t$ if and only if $C(s)\subseteq C(t)$. Let us show that $A$ indeed is an atomization for $F_{C;T}(\emptyset)$.

    Axioms \eqref{axiom:AS1}, \eqref{axiom:AS2} and \eqref{axiom:AS4} and \eqref{axiom:AS5} all either follow immediately or are assumptions on the atomization itself. To show that \eqref{axiom:AS6} holds, let $t$ be a term. Then, there always exists $c\in C(t)$, and $\psi_c < c \le t$. It remains to prove that \eqref{axiom:AS3} characterizes the order relation of $F_{C;T}(\emptyset)$.

    Assume that $(s,t)$ is a pair of terms such that $\neg\exists\phi\colon \big((\phi<s)\wedge (\phi\nless t)\big)$. Equivalently, for every $c\in C$, if $\psi_c < s$, then $\psi_c < t$. This is equivalent to saying that if $c\in C(s)$, then $c\in C(t)$. But this means that $C(s)\subseteq C(t)$, thus $F_{C;T}(\emptyset)\models s\le t$.

    Assume now that $F_{C;T}(\emptyset)\models s\le t$, thus $C(s)\subseteq C(t)$. Then, if $\psi_c < s$, $c\in C(s)\subseteq C(t)$, and therefore $\psi_c < t$. This shows that $\neg\exists\phi\colon \big((\phi<s)\wedge (\phi\nless t)\big)$, and thus we conclude that \eqref{axiom:AS3} holds.
\end{proof}

\begin{proposition}\label{proposition:atoms_and_freest_model}
    Let $M$ be a semilattice over a set of constants $C$ with terms $T$, regarded as atomized by $\Omega(M)$, which atomizes every semilattice, Theorem \ref{theorem:every_semilattice_atomizable}. Let $t, s\in F_{C;T}(\emptyset)$ be two terms. Let $\phi$ be an atom over $C$ and $\nu_M\colon F_{C;T}(\emptyset)\to M$ the natural homomorphism.
    \begin{enumerate}
        \item $\big(F_{C;T}(\emptyset)\models (t\le s)\big) \Rightarrow M\models \big(\nu_M(t)\le \nu_M(s)\big)$;
        \item $(\phi\in M)\wedge \big(F_{C;T}(\emptyset)\models (\phi<s)\big)\Leftrightarrow M\models \big(\phi<\nu_M(s)\big)$.
    \end{enumerate}
\end{proposition}

\begin{proof}
    (1) is an immediate consequence of $\nu_M$ being a homomorphism. (2) follows from Proposition \ref{proposition:atomizations_and_order} (3), since the atoms of $M$ are exactly the atoms of $\Omega(M)$.
\end{proof}

The following result is now immediate from Proposition \ref{proposition:atoms_and_freest_model}.

\begin{proposition}
    Let $M$ and $N$ be two semilattices over a set of constants $C$ with respective terms $T$ and $T'$ and let $t, s\in T\cap T'$. If an atom $\phi$ discriminates a duple $(t, s)$ in $M$, then either $\phi\not\in N$ or $\phi$ discriminates $(t, s)$ in $N$.
\end{proposition}

\begin{proof}
    The result follows from applying Proposition \ref{proposition:atoms_and_freest_model} (2) to $M$ and to $N$.
\end{proof}

One important atom is the zero atom, $\ominus_C$, defined as the atom whose upper constant segment is $C$. This atom discriminates no duple, so it is compatible with every semilattice over $C$. We say that a set of atoms \emph{spawns} a semilattice $M$ with terms $T$ if it atomizes $M$; the two words are synonyms, and every condition stated for one applies to the other. Furthermore, we have the following.

\begin{proposition}\label{proposition:spawned}
    Let $A$ be a set of atoms over the set of constants $C$ and let $T$ be a set of terms over $C$. Then $A\cup\{\ominus_C\}$ atomizes a semilattice $M_A$ over $C$ with terms $T$, unique up to an interpretation-preserving isomorphism, Definition \ref{definition:interpretation_preserving}. Moreover, $A$ itself atomizes $M_A$ if and only if the upper constant segments of the atoms of $A$ cover $C$.
\end{proposition}

\begin{proof}
    For $t\in T$ let $L_A(t)=\{\phi\in A\colon U^c(\phi)\cap C(t)\ne\emptyset\}$. Since $C(t\cup s) = C(t)\cup C(s)$,
    \begin{equation*}
        L_A(t\cup s) = L_A(t)\cup L_A(s) \qquad\text{and}\qquad L_A(t) = \bigcup_{c\in C(t)} L_A(\{c\}).
    \end{equation*}
    Let $M_A = \{L_A(t)\colon t\in T\}$. Since $T$ is closed under finite unions, the first identity shows that $M_A$ is closed under union, so it is a semilattice with union as join and inclusion as order, and $\nu(t) = L_A(t)$ is a surjective homomorphism $F_{C;T}(\emptyset)\to M_A$. By the second identity, $\nu(t)$ contains $\nu(\{c\})$ for every $c\in C(t)$, and every element of $M_A$ containing all of them contains $\nu(t)$. That is, $\nu(t)$ is the least upper bound in $M_A$ of the constants of $t$, also when $t$ is infinite. Hence $M_A$ is a semilattice over $C$, with the constant $c$ interpreted as $\nu(\{c\})$, and $T$ is a set of terms for $M_A$, Definition \ref{definition:admissible_terms}.

    For $\phi\in A\cup\{\ominus_C\}$, set $\phi<\nu(t)$ if and only if $\phi<t$, as in Definition \ref{definition:atom_from_ucs}. This depends only on the element $\nu(t)$, since for $\phi\in A$ it says $\phi\in L_A(t)$, and $\ominus_C$ is lower than every element. It also extends the order of $M_A$: if $\phi\in A$, $\phi<\nu(t)$ and $\nu(t)\le\nu(s)$, then $\phi\in L_A(t)\subseteq L_A(s)$, so $\phi<\nu(s)$. Then \eqref{axiom:AS1} holds because upper constant segments are non-empty, \eqref{axiom:AS2} because atoms are only related to regular elements from below, \eqref{axiom:AS4} because $U^c(\phi)$ meets $C(t)$ if and only if it contains some $c\in C(t)$, and \eqref{axiom:AS5} because an atom is determined by its upper constant segment. For \eqref{axiom:AS3}, $\nu(t)\not\le\nu(s)$ if and only if some $\phi\in A$ is in $L_A(t)\setminus L_A(s)$, that is, discriminates $(t,s)$, while $\ominus_C$ discriminates no duple. Finally, \eqref{axiom:AS6} holds for $A\cup\{\ominus_C\}$ because $\ominus_C$ is lower than every element, and it holds for $A$ if and only if every constant is above some atom of $A$, since every element is $\nu(t)$ for some $t\in T$ and $t$ contains a constant.

    As for uniqueness, let $N$ be a semilattice over $C$ with terms $T$ atomized by $A\cup\{\ominus_C\}$, or by $A$. By Proposition \ref{proposition:atomizations_and_order} (3) and (5), $\nu_N(t)\le\nu_N(s)$ if and only if $L_A(t)\subseteq L_A(s)$, the zero atom being lower than both or than neither. Hence $\nu_N(t)\mapsto L_A(t)$ is a well-defined order isomorphism $N\to M_A$, and it is interpretation-preserving.
\end{proof}

    We have thus shown that every set of atoms that contains $\ominus_C$, or whose upper constant segments cover $C$, spawns a semilattice with terms $T$, unique up to isomorphism. We have chosen to use the elements of $F_{C;T}(\emptyset)$ to name and keep track of elements of a semilattice $M$ since it allows for this framework of atomizations. In this way, atoms can be carried to different models, allowing us to transfer properties across them.

Following \cite{MarPol21}*{Theorem 37}, atomizations correspond to subdirect product decompositions into two-element semilattices. Recall, \cite{BurSan81}*{Definition 8.1}, that a semilattice is a \emph{subdirect product} of a family $(M_i)_{i\in I}$ if it is a subalgebra of $\prod_{i\in I}M_i$ that projects onto every factor.

\begin{theorem}\label{theorem:subdirect}
    Let $M$ be a semilattice over $C$ with terms $T$, and assume that $M$ has at least two elements. For an atom $\phi\ne\ominus_C$ over $C$, let $M_\phi$ be the semilattice spawned by $\{\phi,\ominus_C\}$ with terms $T$ through Proposition \ref{proposition:spawned}. In particular, $M_\phi$ has two elements $0<1$, with $\phi<1$ and $\phi\nless0$.
    \begin{enumerate}
        \item If $A$ atomizes $M$ and $A'=A\setminus\{\ominus_C\}$, then $A'\ne\emptyset$. For $\phi\in A'$ the map $\chi_\phi\colon M\to M_\phi$, $\chi_\phi(x)=1$ if $\phi<x$ and $\chi_\phi(x)=0$ otherwise, is an interpretation-preserving homomorphism and $x\mapsto(\chi_\phi(x))_{\phi\in A'}$ is an isomorphism of $M$ onto a subdirect product of $(M_\phi)_{\phi\in A'}$.
        \item Conversely, let $(M_i)_{i\in I}$ be two-element semilattices $\{0<1\}$ over $C$ and $\pi_i\colon M\to M_i$ interpretation-preserving homomorphisms such that $x\mapsto(\pi_i(x))_{i\in I}$ is injective. Set $U_i=\{c\in C\colon\nu_{M_i}(c)=1\}$ and $\phi_i:<U_i$. Then $\phi_i\ne\ominus_C$, the set $\{\phi_i\colon i\in I\}\cup\{\ominus_C\}$ atomizes $M$, with $\phi_i<x$ if and only if $\pi_i(x)=1$, and $M_i=M_{\phi_i}$, $\pi_i=\chi_{\phi_i}$.
    \end{enumerate}
\end{theorem}

\begin{proof}
    By Proposition \ref{proposition:atomizations_and_order} (3), for $t\in T$ the segment $L^a_{M_\phi}(\nu_{M_\phi}(t))$ is $\{\phi,\ominus_C\}$ or $\{\ominus_C\}$. Both occur, the first for $t=\{c\}$ with $c\in U^c(\phi)$ and the second for $t=\{c'\}$ with $c'\in C\setminus U^c(\phi)$, which is non-empty because $\phi\ne\ominus_C$. By Proposition \ref{proposition:atomizations_and_order} (5) the lower atomic segments determine the elements as well as the order on them. Hence $M_\phi$ has exactly two elements $0<1$, with $\phi<1$ and $\phi\nless0$. Moreover, $\nu_{M_\phi}(t)=1$ if and only if $C(t)\cap U^c(\phi)\ne\emptyset$.

    (1) Some duple of $M$ is negative, so by \eqref{axiom:AS3} some atom of $A$ discriminates the duple. That atom is not $\ominus_C$, which is lower than every regular element, so $A'\ne\emptyset$. For $\phi\in A'$ and $t\in T$, Proposition \ref{proposition:atomizations_and_order} (3) gives $\phi<\nu_M(t)$ in $M$ if and only if $C(t)\cap U^c(\phi)\ne\emptyset$, if and only if $\nu_{M_\phi}(t)=1$. In other words, $\chi_\phi(\nu_M(t))=\nu_{M_\phi}(t)$. Since both $\nu_M$ and $\nu_{M_\phi}$ are surjective homomorphisms, this identity says that $\chi_\phi$ is a surjective interpretation-preserving homomorphism. Hence $x\mapsto(\chi_\phi(x))_{\phi\in A'}$ is a homomorphism into $\prod_{\phi\in A'}M_\phi$, the operation being componentwise, and its composite with the projection onto $M_\phi$ is the surjection $\chi_\phi$. So the image of this map is a subdirect product of the $M_\phi$. The map is also injective. If $\chi_\phi(x)=\chi_\phi(y)$ for every $\phi\in A'$, then $\phi<x$ if and only if $\phi < y$. In other words, $L^a_M(x)=L^a_M(y)$, since $\ominus_C$ is in both, yielding $x=y$.

    (2) Each element of $M_i$ is a join of constants, so $\nu_{M_i}$ takes both values on $C$ and $\emptyset\ne U_i\ne C$. Hence $\phi_i$ is an atom over $C$ and $\phi_i\ne\ominus_C$. Since $\{0<1\}$ is complete, $T$ is a set of terms for $M_i$, and $\nu_{M_i}(t)=1$ if and only if $C(t)\cap U_i\ne\emptyset$. So $M_i=M_{\phi_i}$. As $\pi_i$ preserves the interpretation of terms, $\pi_i(\nu_M(t))=1$ if and only if $C(t)\cap U_i\ne\emptyset$.

    We observe next that, for $x,y\in M$, $x\le y$ if and only if $\pi_i(x)\le\pi_i(y)$ for every $i\in I$. Indeed, homomorphisms are order-preserving. Conversely, if $\pi_i(x)\le\pi_i(y)$ for every $i$, then $\pi_i(x\odot y)=\pi_i(x)\odot\pi_i(y)=\pi_i(y)$ for every $i$, and $x\odot y=y$ because $x\mapsto(\pi_i(x))_{i\in I}$ is injective.

    Now let $N$ be the semilattice spawned by $B=\{\phi_i\colon i\in I\}\cup\{\ominus_C\}$ with terms $T$. For $t\in T$, by Proposition \ref{proposition:atomizations_and_order} (3), $L^a_N(\nu_N(t))$ consists of $\ominus_C$ and of the $\phi_i$ with $C(t)\cap U_i\ne\emptyset$, that is, with $\pi_i(\nu_M(t))=1$. Hence, by Proposition \ref{proposition:atomizations_and_order} (5), $\nu_N(t)\le\nu_N(s)$ if and only if every $\phi_i$ lower than $\nu_N(t)$ is lower than $\nu_N(s)$, that is, $\pi_i(\nu_M(t))\le\pi_i(\nu_M(s))$ for every $i$. By the previous paragraph, this is $\nu_M(t)\le\nu_M(s)$. Consequently $f(\nu_N(t))=\nu_M(t)$ defines a bijection $f\colon N\to M$ that preserves and reflects the order, that is, an interpretation-preserving isomorphism. Transporting the atomization $B$ of $N$ along $f$, $B$ atomizes $M$, with $\phi_i<\nu_M(t)$ if and only if $\phi_i<\nu_N(t)$, if and only if $\pi_i(\nu_M(t))=1$. So $\phi_i<x$ if and only if $\pi_i(x)=1$ for $x\in M$, and $\chi_{\phi_i}=\pi_i$.
\end{proof}

Thus each atomization of $M$ exhibits $M$ as a subdirect product of two-element semilattices, one factor per atom other than $\ominus_C$, and each factor map $\chi_\phi$ preserves the joins named by $T$, not only the finite ones. In the finite case this is \cite{MarPol21}*{Theorem 37}. Conversely, every subdirect representation of $M$ by two-element semilattices over $C$ through interpretation-preserving maps arises in this way. Moreover, the two constructions are mutually inverse up to $\ominus_C$ and up to repeated factors. Indeed, (2) applied to $(\chi_\phi)_{\phi\in A'}$ returns $A'\cup\{\ominus_C\}$, and (1) applied to the atomization of (2) returns $(\pi_i)_{i\in I}$ with equal maps identified.

The requirement that the $\pi_i$ preserve the interpretation of the terms of $T$ cannot be weakened to their being homomorphisms. In the ascending chain $c_0<c_1<\cdots<c_\infty$, with $M$ consisting only of the constants and all non-empty subsets of $C$ as terms, the map sending every $c_i$ to $0$ and $c_\infty$ to $1$ is a surjective homomorphism, but the atom $\phi_\infty:<\{c_\infty\}$ it would define is not compatible with $M$, since $\nu_M(t)=c_\infty$ for $C(t)=\{c_i\mid i\in\mathbb{N}\}$.

\section{Freedom and full crossing}\label{section:full_crossing}

In this section we introduce the join of atoms, Definition \ref{definition:atom_join}, and the concept of relative freedom of models, Definition \ref{definition:freer_or_as_free}. We then focus on the introduction of the full crossing operation. We first introduce the full crossing of one duple in a model, Definition \ref{definition:full_crossing}, which is a direct generalization of full crossing in the case of finite atomized semilattices, \cite{MarPol21}*{Definition 3.2}. Just like in the finite case, we then show that full crossing produces the freest semilattice that models the positive theory of the initial semilattice together with the duple being crossed, Theorem \ref{theorem:full_crossing_freest}. Equivalently, the result of crossing a duple $r = (r_L, r_R)$ in a model $M$ is the quotient of $M$ by the principal congruence $\Theta(r_R, r_L \odot r_R)$, Theorem \ref{theorem:congruences}.

The procedure introduced in Definition \ref{definition:full_crossing} performs atom joins and removals when crossing each duple, in such a way that a limit atomization may be difficult to obtain when crossing an infinite set of duples. We will nonetheless show that a limit model does exist. Indeed, one way of obtaining such limit is using atomizations with all the compatible atoms (full atomizations), which lets us introduce an equivalent procedure for full atomizations, Definition \ref{definition:full_crossing_infinite}. We are then able to show that every model is atomizable, Theorem \ref{theorem:every_semilattice_atomizable}, and we can then conclude that the crossing procedure works for infinite sets of duples, Proposition \ref{proposition:full_crossing_freest}.

We finish this section by showing that the full crossing procedure is equivalent to computing a quotient of the given algebra by the least $T$-aligned congruence containing the duples, Theorem \ref{theorem:congruences}.

\begin{definition}
    Let $\phi$ and $\psi$ be two atoms. We say that $\phi$ is \emph{wider} than $\psi$, denoted $\phi < \psi$, if $U^c(\phi)\supsetneq U^c(\psi)$.
\end{definition}

\begin{definition}\label{definition:atom_join}
    Let $\phi$ and $\psi$ be two atoms. The \emph{join} of $\phi$ and $\psi$, denoted by $\phi\join\psi$, is the atom with upper constant segment $U^c(\phi\join\psi) = U^c(\phi)\cup U^c(\psi)$. More generally, if $\phi_i$, $i \in I$, are atoms with $I\ne\emptyset$, $\join_{i\in I} \phi_i$ is the atom with upper constant segment $U^c(\join_{i\in I} \phi_i) = \cup_{i\in I} U^c(\phi_i)$.
\end{definition}

The join of atoms has the following properties.

\begin{proposition}\label{proposition:join_properties}
    Let $x$ be a term over a set of constants $C$ and let $\phi$, $\psi$ and $\phi_i$, $i\in I$, $I\ne\emptyset$, be atoms over $C$, with $<$ read as in Definition \ref{definition:atom_from_ucs}. Then,
    \begin{enumerate}
        \item $\join$ is idempotent, commutative and associative;
        \item $\phi < x \Rightarrow (\phi\join\psi) < x$;
        \item $(\join_{i\in I} \phi_i < x) \Leftrightarrow(\exists i\in I \colon \phi_i < x)$;
        \item $(\phi\join\psi < x)\wedge (\phi\nless x)\Rightarrow (\psi < x)$;
        \item $\phi$ is wider or equal to $\psi$ if and only if $\phi = \phi\join\psi$.
    \end{enumerate}
\end{proposition}

\begin{proof}
    By Definition \ref{definition:atom_from_ucs}, $\phi < x$ if and only if $U^c(\phi)\cap C(x)\ne\emptyset$, and by Definition \ref{definition:atom_join}, $U^c(\join_i\phi_i) = \cup_i U^c(\phi_i)$. Items (2)--(5) are immediate from these two facts, and (1) follows from the corresponding properties of the union, since an atom is determined by its upper constant segment.
\end{proof}

If $x = \nu_M(t)$ and the atoms belong to an atomization of a semilattice $M$, this is the relation $<$ of $M$, as shown after Definition \ref{definition:atom_from_ucs}. Moreover, joining atoms compatible with $M$ always produces another atom compatible with $M$, as the next corollary shows.

\begin{corollary}\label{corollary:join_still_in_model}
    Let $M$ be an atomized semilattice over $C$ and let $\{\psi_i\}_{i\in I}$, $I\ne\emptyset$, be a set of atoms compatible with $M$. Then $\join_i\psi_i$ is compatible with $M$.
\end{corollary}

\begin{proof}
    Assume that $r = (r_L, r_R)$ is a duple of $M$ discriminated by $\join_i\psi_i$. Then $\join_i\psi_i\nless r_R$, from which $U^c(\join_i\psi_i)\cap C(r_R) = \emptyset$. Since $U^c(\join_i\psi_i) = \cup_{i\in I} U^c(\psi_i)$, $U^c(\psi_i)\cap C(r_R) = \emptyset$, $\forall i$. On the other hand, $\join_i\psi_i < r_L$, hence by Proposition \ref{proposition:join_properties} (3), $\exists i\colon \psi_i < r_L$. We then deduce that $\psi_i$ discriminates $r$, thus $M\models r^-$. Therefore, $\join_i\psi_i$ is compatible with $M$.
\end{proof}

We now move on to introducing the concept of relative freedom for a pair of semilattices, an idea that will lead us to introduce the full crossing operator.

\begin{definition}\label{definition:freer_or_as_free}
    Let $M$ and $N$ be two semilattices over a set of constants $C$ with terms $T$. We say that $M$ is \emph{freer or as free} as $N$ if there is a surjective, interpretation-preserving homomorphism $f\colon M\to N$, Definition~\ref{definition:interpretation_preserving}.
\end{definition}

For example, the freest finitely generated semilattice $F_C(\emptyset)$ is freer or as free as any other finitely generated semilattice over $C$. Relative freedom can be characterized in a simple way, in correspondence to how it was defined in \cite{MarPol21}.

\begin{proposition}\label{proposition:characterization_of_freedom}
    Let $M$ and $N$ be two semilattices over $C$ with equal sets of terms $T$. Then, $M$ is freer or as free as $N$ if and only if for every duple $r$ of $M$, $N\models r^-\Rightarrow M\models r^-$.
\end{proposition}

\begin{proof}
    Assume first that $M$ is freer or as free as $N$, and let $f\colon M\to N$ be a surjective, interpretation-preserving homomorphism. Let $r = (r_L, r_R)$ be a duple of $M$ with $M\models r^+$, that is, $\nu_M(r_L)\le\nu_M(r_R)$. Since $f$ is order-preserving and $f(\nu_M(t)) = \nu_N(t)$ for $t\in T$, we get $\nu_N(r_L)\le\nu_N(r_R)$, that is, $N\models r^+$. Thus $N\models r^-$ implies $M\models r^-$.

    Assume now that $N\models r^-\Rightarrow M\models r^-$. This is equivalent to $M\models r^+\Rightarrow N\models r^+$. Consider the map $f\colon M\to N$ sending the element of $M$ represented by a term $t\in T$ to the element of $N$ represented by the same term $t$, that is, $f(\nu_M(t)) = \nu_N(t)$. This map is well defined, since $M\models (t = s)$ implies $N\models (t = s)$ by hypothesis. It is surjective, and it is interpretation-preserving by construction. Using $t\cup s\in T$ and that $\nu_M$ and $\nu_N$ are homomorphisms, we get $f(\nu_M(t)\odot\nu_M(s)) = f(\nu_M(t\cup s)) = \nu_N(t\cup s) = \nu_N(t)\odot\nu_N(s)$, thus $f$ is a homomorphism.
\end{proof}

This is thus compatible with \cite{MarPol21}*{Definition 3.1}. In particular, these homomorphisms are the ones obtained in \cite{MarPol21}*{Theorem 28 i)}. We define relative freedom in terms of homomorphisms here since, unlike in the finite case, different infinite semilattices over a given set of constants may not admit the same joins of constants.

We now introduce some terminology.

\begin{definition}
    Let $M$ be an atomized semilattice over $C$. The \emph{discriminant} of a duple $r = (r_L, r_R)$ of $M$ is the set $\dis_M(r) = L_M^a(r_L)\setminus L_M^a(r_R)$.  If $R$ is a set of duples of $M$, the discriminant of $R$ is the set $\dis_M(R) = \cup_{r\in R}\dis_M(r)$.
\end{definition}

Notice that, as a consequence of \eqref{axiom:AS3}, $M\models r^+$ if and only if $\dis_M(r) = \emptyset$.

\begin{definition}
    Let $M$ be a semilattice over $C$. The \emph{theory of $M$}, $\Th_0(M)$, is the set of duples of $M$ satisfied by $M$, that is, $\Th_0(M) = \Th_0^+(M)\cup\Th_0^-(M)$, where the \emph{positive and negative theories of $M$}, $\Th_0^+(M)$ and $\Th_0^-(M)$, are respectively the sets of positive and negative duples satisfied by $M$.
\end{definition}

Given $M$ and $N$ two atomized semilattices over a set of constants $C$ with the same terms $T$ and respectively atomized by sets $A$ and $B$, we can produce a new model $M+N$ over the same terms and atomized by $A\cup B$, the semilattice spawned by $A\cup B$, Proposition \ref{proposition:spawned}, since the upper constant segments of the atoms of $A$ already cover $C$. We then have the following.

\begin{proposition}
    Let $M$ and $N$ be two models over $C$ with terms $T$ and atomized by respective sets of atoms $A$ and $B$. Then, $\Th_0^+(M+N) = \Th_0^+(M)\cap\Th_0^+(N)$ and $\Th_0^-(M+N) = \Th_0^-(M)\cup \Th_0^-(N)$.
\end{proposition}

\begin{proof}
    Let $r$ be a duple of $M$, equally of $N$ or of $M+N$, since the three have terms $T$. By \eqref{axiom:AS3}, $M+N\models r^-$ if and only if some atom of $A\cup B$ discriminates $r$. By Proposition \ref{proposition:atoms_and_freest_model}~(2), whether an atom discriminates $r$ does not depend on the semilattice in which $r$ is considered. Hence, by \eqref{axiom:AS3} for $M$ and for $N$, $M+N\models r^-$ if and only if $M\models r^-$ or $N\models r^-$. This is $\Th_0^-(M+N) = \Th_0^-(M)\cup\Th_0^-(N)$, and taking complements in the set of duples of $T$ gives $\Th_0^+(M+N) = \Th_0^+(M)\cap\Th_0^+(N)$. This is the argument of \cite{MarPol21}*{Theorem 14}.
\end{proof}

One important question to ask is the following. If we have a semilattice in which a certain duple is negative, how can we make that duple positive while changing the least possible amount of duples? With that purpose, we introduce the following operation.

\begin{definition}\label{definition:full_crossing}
    Let $M$ be an atomized semilattice over a set of constants $C$ and with terms $T$, and let $r = (r_L, r_R)$ be a duple of $M$. Assume that $M$ is atomized by $A$. The \emph{full crossing of $r$ in $M$}, denoted $\crossing_r M$, is the model with the same terms as $M$ and atomized by $(A\setminus H)\cup (H\join B)$, where $H = \dis_M(r)$, $B = L_M^a(r_R)$ and $H\join B = \{\phi\join\psi\colon (\phi\in H)\wedge(\psi\in B)\}$.
\end{definition}

If $r = (r_L, r_R)$, we also write $\crossing_r M = \crossing_{r_L\le r_R}M$. Notice that $(A\setminus H)\cup (H\join B)$ is an atomization for a semilattice with the same term set as $M$, by Proposition \ref{proposition:spawned}, provided that the upper constant segments of its atoms cover $C$, that is, that it satisfies \eqref{axiom:AS6}. It does, because $A$ also does. Indeed, if we have $a\in M$, there is an atom $\phi\in A$ such that $\phi < a$. If this atom is removed, an atom of the form $\phi\join\psi$ is in $\crossing_r M$, with $\psi\in L_M^a(r_R)\ne\emptyset$. And by Proposition \ref{proposition:join_properties} (2), $\phi\join\psi < a$. Furthermore, none of the atoms in $H\join B$ discriminate $r$. Since every atom in $\dis_M(r)$ has been removed, $\crossing_r M \models r^+$. In fact, we can compute exactly what duples become positive when performing full crossing.

\begin{proposition}\label{proposition:logical_cosequences}
    Let $M$ be an atomized semilattice over a set of constants $C$ with terms $T$ and let $r=(r_L, r_R)$ and  $s=(s_L, s_R)$ be duples of $M$. If $M\models s^-$ and $\crossing_r M\models s^+$, then
    \[M\models \big((s_L \le s_R\odot r_L)\wedge (r_R\le s_R)\big).\]
\end{proposition}

\begin{proof}
    Assume that $M$ is atomized by a set of atoms $A$. First notice that, since $M\models s^-$ and $\crossing_r M\models s^+$, every atom of $A$ discriminating $s$ must be removed when performing full crossing. Therefore, $\dis_M(s)\subseteq\dis_M(r)$. This implies that every atom discriminating $s$ is below $r_L$, from which we deduce that $M\models s_L\le s_R\odot r_L$.

    Now take $H=\dis_M(r)$ and $B = L_M^a(r_R)$. Since $M\models s^-$, there is at least one atom $\phi\in\dis_M(s)$, which also discriminates $r$. For every $\psi\in B$, the atom $\phi\join\psi$ is in $\crossing_r M$, and by Proposition \ref{proposition:join_properties} (3), $\phi\join\psi$ is in $s_L$. Since $\crossing_r M\models s^+$, $\phi\join\psi$ is in $s_R$. By Proposition \ref{proposition:join_properties} (4), since $\phi\nless s_R$ and $\phi\join\psi < s_R$, $\psi < s_R$, for every $\psi\in B$. But $B = L_M^a(r_R)$, thus $L_M^a(r_R)\subseteq L_M^a(s_R)$, and then $M\models r_R\le s_R$.
\end{proof}

We now introduce the main result involving full crossing. But first, let us introduce the following terminology.

\begin{definition}\label{definition:freest_over_duples}
    Let $T$ be a set of terms over a set of constants $C$. Let $R$ be a set of duples of elements in $T$. Then, $F_{C;T}(R)$ denotes the freest semilattice over $C$ with terms $T$ that models $R^+$. Similarly, $F_C(R)$ and $\tilde{F}_C(R)$ respectively denote the freest finitely-generated and complete semilattices that model $R^+$.
\end{definition}

Consequently, the positive theory of $F_{C;T}(R)$ exclusively contains $R$ and the logical consequences of $R$. Existence and uniqueness of $F_{C;T}(R)$ are shown in Corollary \ref{corollary:freest_models}.

\begin{theorem}\label{theorem:full_crossing_freest}
    Let $M$ be an atomized semilattice over a set of constants $C$ with terms $T$ and let $r$ be a duple of $M$. The full crossing of $r$ in $M$ is the freest model of $\Th_0^+(M)\cup r^+$ with terms $T$, $\crossing_r M = F_{C;T}(\Th_0^+(M)\cup \{r\})$.
\end{theorem}

\begin{proof}
    The proof is very similar to that of the finite case, \cite{MarPol21}*{Theorem 16}. However, we will still introduce it here as it is quite illustrative. Since we are only changing the atomization, the admissible terms remain the same. We thus have to prove that the resulting atomization provides a semilattice modeling $\Th_0^+(M)\cup r^+$ and that every duple positive in it is positive in every semilattice with terms $T$ modeling $\Th_0^+(M)\cup r^+$.

    Let $r = (r_L, r_R)$. Let $A$ be the atomization of $M$ and let $H = \dis_M(r)$. Let $B = L_M^a(r_R)$ be the atoms in $A$ that are below $r_R$. The atomization of $\crossing_r M$ is then $(A\setminus H)\cup (H\join B)$, where $H\join B = \{\phi\join\psi\colon (\phi\in H)\wedge(\psi\in B)\}$. As mentioned above, $\crossing_r M \models r^+$.

    By Corollary \ref{corollary:join_still_in_model}, the join of two atoms compatible with a semilattice is an atom compatible with that semilattice. Thus, $H\join B\subseteq \Omega(M)$, and thus they do not discriminate any duple which is positive in $M$. Therefore, $\crossing_r M\models \Th_0^+(M)$.

    Let $N$ be a semilattice with terms $T$ such that $N\models\Th_0^+(M)\cup r^+$. A duple positive in $\crossing_r M$ and in $M$ is positive in $N$. It remains to prove that if $s$ is a duple of $M$ such that $M\models s^-$ and $\crossing_r M\models s^+$, then $N\models s^+$. Let $s$ be one such duple. By Proposition \ref{proposition:logical_cosequences} we know that
    \[M\models \big((s_L \le s_R\odot r_L)\wedge (r_R\le s_R)\big).\]
    Since $T$ is closed under finite unions, $(s_L, s_R\cup r_L)$ and $(r_R, s_R)$ are duples with terms in $T$, so they belong to $\Th_0^+(M)$ and hold in $N$. Then, in $N$, $s_L \le s_R \odot r_L \le s_R \odot r_R = s_R$.
    Therefore, $N\models s^+$, and by Proposition \ref{proposition:characterization_of_freedom}, $\crossing_r M$ is freer or as free as $N$.
\end{proof}

Notice that no assumption is made on the duple $r$; if $M\models r^+$, then $H=\emptyset$ and $\crossing_r M = M$. We can now easily prove the following.

\begin{corollary}\label{corollary:cardinal_atomization_crossing}
    Let $M$ be an atomized semilattice over a constant set $C$ with terms $T$, atomized by an atomization $A$ with $|A|=\kappa$ an infinite cardinal. Let $r$ be any duple of $M$. Then $F_{C;T}(\Th_0^+(M)\cup\{r\})$ admits an atomization of cardinality at most $\kappa$.
\end{corollary}

\begin{proof}
    As a consequence of Theorem \ref{theorem:full_crossing_freest}, $F_{C;T}(\Th_0^+(M)\cup\{r\}) = \crossing_r M$. By definition, $\crossing_r M$ is atomized by $(A\setminus\dis_M(r))\cup (H\join B)$, where $H = \dis_M(r)$ and $B = L_M^a(r_R)$ are subsets of $A$. Now, since $H\join B$ has at most one different atom per each pair of one atom in $H$ and another in $B$, $|H\join B| \le |H||B| \le \kappa\cdot\kappa = \kappa$. The considered atomization for the full crossing then has cardinality $|A\setminus\dis_M(r)| + |H\join B|\le \kappa + \kappa = \kappa$.
\end{proof}

Now we raise the question on how we can use the full crossing procedure to obtain the freest model satisfying the theory of a given model as well as a set $R$ of duples. First, if $R$ is finite, we can just apply the procedure as introduced in Definition \ref{definition:full_crossing} a finite number of times. We then have to wonder if the resulting semilattice is independent of the order in which we cross the duples. And indeed, using Theorem \ref{theorem:full_crossing_freest} we can prove the following.

\begin{corollary}\label{corollary:commutativity}
    The full crossing operation over a finite set of duples is commutative. More concretely, although the resulting atomizations may be different, the resulting semilattice is independent of the order of crossings.
\end{corollary}

\begin{proof}
    Assume we have an atomized semilattice $M$ over a set of constants $C$ with terms $T$, and let $R$ be a finite set of duples of $M$. Assume we have fixed an order $R=\{r_1,r_2,\dots,r_n\}$. A semilattice with terms $T$ models $\Th_0^+(F_{C;T}(X))$ if and only if it models $X$, for any set of duples $X$, by Proposition \ref{proposition:characterization_of_freedom}, so crossing $r_{k+1}$ in $F_{C;T}(\Th_0^+(M)\cup\{r_1,\dots,r_k\})$ gives $F_{C;T}(\Th_0^+(M)\cup\{r_1,\dots,r_{k+1}\})$, by Theorem \ref{theorem:full_crossing_freest} and the uniqueness of freest models, Corollary \ref{corollary:freest_models} (1). By a repeated application of this fact we deduce that $\crossing_{r_n}\dots\crossing_{r_2}\crossing_{r_1} M = F_{C;T}(\Th_0^+(M)\cup R)$, and thus, the resulting semilattice is independent of the order chosen for the duples of $R$.
\end{proof}

We can thus introduce the following.

\begin{definition}\label{definition:full_crossing_finite}
    Let $M$ be an atomized semilattice over a set of constants $C$ with terms $T$ and let $R = \{r_1,r_2,\dots,r_n\}$ be a finite set of duples of $M$. We denote $\crossing_R M = \crossing_{r_n}\dots\crossing_{r_2}\crossing_{r_1} M$.
\end{definition}

Now, since in the full crossing procedure for a set of duples $R$ we are iteratively removing and creating new atoms, the limit may not exist if $R$ is infinite. We can get around this by using the full atomization for the model. Indeed, if $M$ is a semilattice over $C$ with terms $T$ atomized by $\Omega(M)$ and $r$ is a duple of $M$, the set $H\join B$ introduced in Definition \ref{definition:full_crossing} consists of joins of atoms of $M$, and it is thus a subset of $\Omega(M)$, Corollary \ref{corollary:join_still_in_model}. Moreover, no atom of $H\join B$ discriminates $r$, since $\phi\join\psi < r_R$ for $\psi\in B$, Proposition \ref{proposition:join_properties} (2), so $H\join B\subseteq\Omega(M)\setminus\dis_M(r)$. Thus, $\crossing_r M$ is atomized by $\Omega(M)\setminus\dis_M(r)$. More generally, even if $R$ is an infinite set of duples, we can introduce the following.

\begin{definition}\label{definition:full_crossing_infinite}
    Let $M$ be an atomized semilattice over a set of constants $C$ and terms $T$, atomized by the full atomization $\Omega(M)$. Let $R$ be a set of duples of $M$. The \emph{full crossing} of $R$ in $M$, denoted $\crossing_R M$, is the model with admissible terms $T$ atomized by $\Omega(M)\setminus\dis_M(R)$, which contains $\ominus_C$, since $\ominus_C$ discriminates no duple, and thus atomizes a semilattice by Proposition \ref{proposition:spawned}.
\end{definition}

If $R$ is finite, Definitions \ref{definition:full_crossing_finite} and \ref{definition:full_crossing_infinite} give the same semilattice, Corollary \ref{corollary:freest_models} (2).

We now focus on proving that every model is atomizable.

\begin{lemma}\label{lemma:compatible_atom_exists}
    Let $C$ be a set of constants and let $M$ be a semilattice over $C$ with set of terms $T$. Let $s$ be a duple of $M$ such that $M\models s^-$. Then, there exists an atom $\phi$ compatible with $M$ such that $\phi$ discriminates $s$.
\end{lemma}

\begin{proof}
    Write $s = (s_L, s_R)$ and consider the set $K = \{c\in C\colon M\models (c\not\le s_R)\}$. If $K$ is empty then every constant $c \in C$ is lower or equal than (the element of $M$ represented by) $s_R$; particularly, every component constant of the term $s_L$ is lower than $s_R$ in $M$, which implies $M\models s^+$ against our assumptions. Since $K$ is not empty, consider the atom $\phi :< K$. We shall prove that this is the desired atom.

    We just argued that if every component constant in $C(s_L)$ is not in $K$ then $M\models s^+$, so there is at least one component constant of $s_L$ in $K$ and then $\phi < s_L$. In addition, since every component constant $c$ of $s_R$ satisfies $c \leq s_R$ in every model, then $c$ is not in $K$, i.e.\ $U^c(\phi)\cap C(s_R) = \emptyset$ or, in other words, $\phi \not< s_R$. Therefore, $\phi$ discriminates $s$.
    We now have to prove that $\phi$ is compatible with $M$. In particular, we will show that if $\phi$ discriminates a duple $r = (r_L, r_R)$, then $M\models r_L\not\le r_R$.

    Assume $\phi$ discriminates $r = (r_L, r_R)$. This requires that $\phi\nless r_R$ which implies $M\models c\le s_R$ for every $c\in C(r_R)$; hence, $M\models r_R\le s_R$. Second, since $\phi < r_L$, we can find $c\in U^c(\phi)\cap C(r_L)$, so $M\models c\not\le s_R$. From this we can deduce that $M\models r_L\not\le s_R$. Suppose $M\models r^+$; we have $M\models r_L\le r_R$ and $M\models r_R\le s_R$ that, by the transitive property implies $M\models r_L\le s_R$, a contradiction.
\end{proof}

\begin{theorem}\label{theorem:every_semilattice_atomizable}
    Let $C$ be a set of constants and let $M$ be a semilattice over $C$ with terms $T$. Then $M$ is atomizable.
\end{theorem}

\begin{proof}
    For every duple $s$ of $M$ such that $M\models s^-$, let $\phi_s$ be an atom discriminating $s$ and compatible with $M$, which exists as a consequence of Lemma \ref{lemma:compatible_atom_exists}. Consider the set $A = \{\phi_s\colon M\models s^-\}\cup\{\ominus_C\}$. Then $A$ is a set of atoms compatible with $M$, the zero atom discriminating no duple at all. Furthermore, $A$ must atomize $M$. Axioms \eqref{axiom:AS1}, \eqref{axiom:AS2}, \eqref{axiom:AS4} and \eqref{axiom:AS5} hold for any set of atoms over $C$ once $\phi<\nu_M(t)$ is read as $U^c(\phi)\cap C(t)\ne\emptyset$, as in Definition \ref{definition:atom_from_ucs}, which depends only on $\nu_M(t)$ because the atoms of $A$ are compatible with $M$, as noted after Definition \ref{definition:compatible}, and \eqref{axiom:AS6} holds because $\ominus_C$ is lower than every regular element. As for \eqref{axiom:AS3}, if $M\models s^-$, there must exist an atom in $A$ discriminating $s$ by construction. On the other hand, since every atom in $A$ is compatible with $M$, if $M\models s^+$, there cannot be any atoms discriminating $s$.
\end{proof}

The zero atom cannot be dispensed with in the proof of Theorem \ref{theorem:every_semilattice_atomizable}. If a constant $c$ is sent to a least element of $M$, then $M\models (c\le s_R)$ for every duple $s = (s_L, s_R)$, so $c\notin K_s$ for every negative $s$, taking $\phi_s :< K_s$ as in the proof of Lemma \ref{lemma:compatible_atom_exists}, where $K_s$ is the set $K$ of that proof, and the family of the $\phi_s$ alone fails \eqref{axiom:AS6} at $\nu_M(c)$.

We can further say the following regarding full crossing.

\begin{proposition}\label{proposition:full_crossing_freest}
    Let $M$ be a semilattice over $C$ with terms $T$ and atomized by $\Omega(M)$ and let $R$ be a set of duples of terms in $T$. Then, $\crossing_R M$ is the freest semilattice model of $\Th_0^+(M)\cup R^+$ with terms $T$.
\end{proposition}

\begin{proof}
    The atoms of $\crossing_R M$ lie in $\Omega(M)$, so none of them discriminates a duple positive in $M$, and none lies in $\dis_M(R)$, so none discriminates a duple of $R$. By \eqref{axiom:AS3}, $\crossing_R M\models(\Th_0^+(M)\cup R^+)$. It remains to show that $\crossing_R M$ is freer or as free as every such model.

    Let $N$ be a semilattice with terms $T$ such that $N\models (\Th_0^+(M)\cup R^+)$. By Theorem \ref{theorem:every_semilattice_atomizable}, we can assume that both $M$ and $N$ are atomized. Let $s$ be a duple such that $N\models s^-$. Then, there exists an atom $\phi\in N$ such that $\phi\in\dis_N(s)$. Since $N\models (\Th_0^+(M)\cup R^+)$, this atom cannot discriminate any duple in $\Th_0^+(M)$, thus we deduce that $\phi\in\Omega(M)$. Furthermore, $\phi$ cannot discriminate any duple in $R$, since $N\models R^+$. Therefore, $\phi\in\Omega(M)\setminus\dis_M(R)$, the atomization of $\crossing_R M$. We deduce that $\crossing_R M \models s^-$. By Proposition \ref{proposition:characterization_of_freedom} this exhibits that $\crossing_R M$ is freer or as free as $N$.
\end{proof}

\begin{corollary}\label{corollary:freest_models}
    Let $T$ be a set of terms over a set of constants $C$.
    \begin{enumerate}
        \item For every set $R$ of duples of $T$, the freest model $F_{C;T}(R)$ of Definition \ref{definition:freest_over_duples} exists and is unique up to an interpretation-preserving isomorphism. It is the semilattice spawned with terms $T$, Proposition \ref{proposition:spawned}, by all the atoms over $C$ that discriminate no duple of $R$, $\ominus_C$ among them.
        \item If $M$ is a semilattice over $C$ with terms $T$ and $R$ is a finite set of duples of $M$, Definitions \ref{definition:full_crossing_finite} and \ref{definition:full_crossing_infinite} give the same semilattice $\crossing_R M$, the freest model of $\Th_0^+(M)\cup R^+$.
    \end{enumerate}
\end{corollary}

\begin{proof}
    (1) Take $M = F_{C;T}(\emptyset)$ in Proposition \ref{proposition:full_crossing_freest}. Every atom over $C$ is compatible with $M$, Proposition \ref{proposition:all_atoms_fit_freest}, so $\Omega(M)$ is the set of all atoms over $C$ and $\crossing_R M$ is spawned by the atoms over $C$ that discriminate no duple of $R$. The positive duples of $M$ are the pairs $(t,s)$ with $C(t)\subseteq C(s)$, which are positive in every semilattice over $C$ with terms $T$, so the models of $\Th_0^+(M)\cup R^+$ are those of $R^+$. Thus $\crossing_R M$ is a freest model of $R^+$. Two freest models of $R^+$ are each freer or as free as the other, so by Proposition \ref{proposition:characterization_of_freedom} they have the same positive duples, and the surjective, interpretation-preserving homomorphisms between them are mutually inverse.

    (2) Both are the freest model of $\Th_0^+(M)\cup R^+$ with terms $T$, by Corollary \ref{corollary:commutativity} and Proposition \ref{proposition:full_crossing_freest} respectively, so they coincide by (1).
\end{proof}

We finish this section by exhibiting that the full crossing can also be seen as a quotient of a semilattice by a principal congruence. Thus, let us first remind the concepts involved.

\begin{definition}
    A congruence on a semilattice $M$, \cite{BurSan81}*{Definition 5.1}, \cite{Pap64, FreNat79}, is an equivalence relation $\theta$ on $M$ that is compatible with the semilattice operation in $M$. More concretely, if $a_1$, $a_2$, $b_1$ and $b_2$ are elements of $M$,
    \[a_1\theta b_1 \wedge a_2\theta b_2\Rightarrow (a_1\odot a_2)\theta(b_1\odot b_2).\]
\end{definition}

If an equivalence relation $\theta$ is a congruence, $M/\theta$ is a semilattice with the operation inherited from $M$. The congruences of a semilattice $M$ form themselves a lattice, where the meet is the intersection of the equivalence relations and the join is the smallest congruence containing the congruences to be joined. Furthermore, it is well known that any congruence can be built from \emph{principal congruences}. If $a,b\in M$, the principal congruence of $a$ and $b$, \cite{BurSan81}*{Definition 5.6}, is the smallest congruence in which $a$ and $b$ are in the same equivalence class. It is denoted by $\Theta(a,b)$. Then, by \cite{BurSan81}*{Theorem 5.7 (d)}, if $\theta$ is any congruence,
\[\theta = \bigcup\big\{\Theta(a,b)\colon (a,b)\in \theta\big\} = \bigvee \big\{\Theta(a,b)\colon (a,b)\in \theta\big\},\]
where $\cup$ is the union of the congruences as sets of pairs and $\vee$ is the join in the lattice of congruences of $M$. We denote by $\pi_\theta\colon M\to M/\theta$ the projection. Since $\pi_\theta$ is a homomorphism, $\pi_\theta(x)\le\pi_\theta(y)$ if and only if $(x\odot y, y)\in\theta$.

When every term of $T$ is finite, as in \cite{MarPol21}, $M/\theta$ is a model of $\Th_0^+(M)$ for every congruence $\theta$, since a homomorphism preserves finite joins. With infinite terms this is no longer the case: $\pi_\theta$ is a homomorphism, and yet $M/\theta$ need not satisfy $\Th_0^+(M)$. The projection $\pi_\theta$ preserves finite joins, but it need not preserve the join $\nu_M(t)$ of the constants of an infinite term $t\in T$. More precisely, let $(r_L, r_R)$ be a positive duple of $M$ and assume $T$ is a set of terms for $M/\theta$. Since $\pi_\theta(\nu_M(r_L))$ is an upper bound of the classes of the constants of $r_L$, we always have $\nu_{M/\theta}(r_L)\le\pi_\theta(\nu_M(r_L))\le\pi_\theta(\nu_M(r_R))$. So the duple can only become negative in $M/\theta$ when $\nu_{M/\theta}(r_R) < \pi_\theta(\nu_M(r_R))$, that is, through an infinite right term. For instance, in the ascending chain $c_0<c_1<\cdots<c_\infty$, with all non-empty subsets of $C$ as terms and $\nu_M(t) = c_\infty$ for the term $t$ with $C(t) = \{c_i\mid i\in\mathbb{N}\}$, let $\theta$ be the congruence with the two classes $\{c_i\mid i\in\mathbb{N}\}$ and $\{c_\infty\}$. Then $M/\theta = \{\bot<\top\}$, the classes of the $c_i$ have join $\bot$, and $\pi_\theta(\nu_M(t)) = \pi_\theta(c_\infty)=\top$. The duple $(c_\infty, t)$ is positive in $M$ but negative in $M/\theta$. The following property of a congruence isolates the congruences for which this does not happen.

\begin{definition}\label{definition:aligned}
    Let $\theta$ be a congruence on a semilattice $M$ over $C$ with terms $T$. Since $c\le\nu_M(t)$ for every $c\in C(t)$, $\pi_\theta(\nu_M(t))$ is an upper bound in $M/\theta$ of the classes $\pi_\theta(c)$, $c\in C(t)$. We say that $\theta$ is \emph{$T$-aligned} if, for every $t\in T$, $\pi_\theta(\nu_M(t))$ is the least upper bound of these classes.
\end{definition}

\begin{lemma}\label{lemma:aligned}
    Let $\theta$ be a congruence on a semilattice $M$ over $C$ with terms $T$. The following are equivalent.
    \begin{enumerate}[label=(\alph*), ref=\alph*]
        \item\label{item:aligned_a} $\theta$ is $T$-aligned.
        \item\label{item:aligned_b} $M/\theta$ is a semilattice over $C$ through $\pi_\theta$, $T$ is a set of terms for $M/\theta$, and $\nu_{M/\theta} = \pi_\theta\circ\nu_M$.
        \item\label{item:aligned_c} $\theta = \ker f$ for a surjective, interpretation-preserving homomorphism $f\colon M\to N$ onto a semilattice $N$ over $C$ with terms $T$.
    \end{enumerate}
    If $\theta$ is $T$-aligned, then $M/\theta\models\Th_0^+(M)$. If every term of $T$ is finite, every congruence on $M$ is $T$-aligned.
\end{lemma}

\begin{proof}
    (\ref{item:aligned_a})$\Rightarrow$(\ref{item:aligned_b}). For $t\in T$, $\pi_\theta(\nu_M(t))$ is the least upper bound of the classes $\pi_\theta(c)$, $c\in C(t)$. As $\nu_M$ and $\pi_\theta$ are surjective, every element of $M/\theta$ is of this form. Hence $T$ is a set of terms for $M/\theta$, Definition \ref{definition:admissible_terms}, with $\nu_{M/\theta}(t) = \pi_\theta(\nu_M(t))$.

    (\ref{item:aligned_b})$\Rightarrow$(\ref{item:aligned_c}). Take $f = \pi_\theta$, which is surjective and, by (\ref{item:aligned_b}), interpretation-preserving.

    (\ref{item:aligned_c})$\Rightarrow$(\ref{item:aligned_a}). By the Homomorphism Theorem, \cite{BurSan81}*{Theorem 6.12}, $f = \bar\iota\circ\pi_\theta$ for an isomorphism $\bar\iota\colon M/\theta\to N$. Hence $\bar\iota(\pi_\theta(\nu_M(u))) = \nu_N(u)$ for every $u\in T$, in particular for $u=t$ and for the singletons $u=\{c\}$. Since $\nu_N(t)$ is the least upper bound of the $\nu_N(c)$, $c\in C(t)$, and $\bar\iota^{-1}$ is an order isomorphism, $\pi_\theta(\nu_M(t))$ is the least upper bound of the $\pi_\theta(c)$, $c\in C(t)$.

    If $\theta$ is $T$-aligned, $\pi_\theta$ shows that $M$ is freer or as free as $M/\theta$, so $M/\theta\models\Th_0^+(M)$ by Proposition \ref{proposition:characterization_of_freedom}. If every term is finite, $\nu_M(t)$ is a finite join of constants, which the homomorphism $\pi_\theta$ preserves.
\end{proof}

The principal congruences generated by two comparable elements are $T$-aligned for an elementary reason: their classes have greatest elements.

\begin{lemma}\label{lemma:principal_aligned}
    Let $M$ be a semilattice over $C$ with terms $T$.
    \begin{enumerate}
        \item If $u\le v$ in $M$, then $x\,\Theta(u,v)\,y$ if and only if $x=y$, or $u\le x$, $u\le y$ and $x\odot v = y\odot v$. In particular, every class of $\Theta(u,v)$ has a greatest element.
        \item Every congruence on $M$ whose classes have greatest elements is $T$-aligned.
    \end{enumerate}
    Consequently, $\Theta(b, a\odot b)$ is $T$-aligned for all $a,b\in M$.
\end{lemma}

\begin{proof}
    (1) Let $\rho$ be the relation described. It is clearly reflexive and symmetric, and it is transitive since $x\odot v=y\odot v=z\odot v$ whenever $x\,\rho\,y\,\rho\,z$ with $x\neq y\neq z$. If $x\,\rho\,y$ and $z\in M$, then $u\le x\odot z$, $u\le y\odot z$ and $x\odot z\odot v = y\odot z\odot v$, so $\rho$ is a congruence, using transitivity and commutativity. It contains $(u,v)$, since $u\le v$ and $u\odot v = v = v\odot v$. Conversely, if $\theta$ is a congruence containing $(u,v)$ and $x\,\rho\,y$ with $x\ne y$, then $x = x\odot u\ \theta\ x\odot v = y\odot v\ \theta\ y\odot u = y$. Hence $\rho = \Theta(u,v)$. If $u\le x$, the class of $x$ contains $x\odot v$, and every $y$ in it satisfies $y\le y\odot v = x\odot v$. Otherwise the class of $x$ is $\{x\}$.

    (2) Let $\theta$ be such a congruence, $t\in T$, and let $\pi_\theta(m)$ be an upper bound of the classes $\pi_\theta(c)$, $c\in C(t)$. Let $\hat m$ be the greatest element of the class of $m$. For $c\in C(t)$, $c\odot\hat m\ \theta\ c\odot m\ \theta\ m\ \theta\ \hat m$, so $c\odot\hat m\le\hat m$, that is, $c\le\hat m$. As $\nu_M(t)$ is the join of the constants of $t$, $\nu_M(t)\le\hat m$, and $\pi_\theta(\nu_M(t))\le\pi_\theta(\hat m)=\pi_\theta(m)$.
\end{proof}

Let $R$ be a set of duples of $M$. Here $\crossing_R M$ is the full crossing of Definition \ref{definition:full_crossing_infinite}, computed from $\Omega(M)$, which atomizes every semilattice; as a semilattice it is the freest model of $\Th_0^+(M)\cup R^+$, Proposition \ref{proposition:full_crossing_freest}, so it depends only on $M$ and $R$ and not on any atomization of $M$. We set $\theta_R = \vee_{(r_L, r_R)\in R}\,\Theta(r_R, r_L\odot r_R)$, where terms stand for their interpretations in $M$. By Proposition \ref{proposition:full_crossing_freest}, $\crossing_R M\models\Th_0^+(M)$, so $M$ is freer or as free as $\crossing_R M$, and Proposition \ref{proposition:characterization_of_freedom} provides a surjective, interpretation-preserving homomorphism $h_R\colon M\to\crossing_R M$. By the Homomorphism Theorem, $\crossing_R M\cong M/\ker h_R$, and we write $\crossing_R M = M/\theta$ when $\theta = \ker h_R$. We have the following result.

\begin{theorem}\label{theorem:congruences}
    Let $M$ be a semilattice over a set of constants $C$ and terms $T$, let $a$ and $b$ be regular elements of $M$ and let $R$ be a set of duples of $M$. Then $\ker h_R$ is the least $T$-aligned congruence on $M$ containing $\theta_R$. In particular, $\crossing_R M = M/\theta_R$ if and only if $\theta_R$ is $T$-aligned. Moreover,
    \begin{enumerate}
        \item $\crossing_{a\le b} M = M / \Theta(b, a\odot b)$;
        \item if $R$ is finite, $\crossing_R M = M/\theta_R$, and every finitely generated congruence on $M$ is $T$-aligned;
        \item $M/\Theta(a,b) = \crossing_{a\le b}\crossing_{b\le a} M$.
    \end{enumerate}
\end{theorem}

\begin{proof}
    For every $(r_L, r_R)\in R$, $\crossing_R M\models (r_R = r_L\odot r_R)$, so $(r_R, r_L\odot r_R)\in\ker h_R$ and $\theta_R\subseteq\ker h_R$. Since $h_R$ is surjective and interpretation-preserving, $\ker h_R$ is $T$-aligned by Lemma \ref{lemma:aligned}. Now let $\theta$ be a $T$-aligned congruence containing $\theta_R$. By Lemma \ref{lemma:aligned}, $M/\theta$ is a semilattice over $C$ with terms $T$, $\nu_{M/\theta} = \pi_\theta\circ\nu_M$ and $M/\theta\models\Th_0^+(M)$. Moreover, $M/\theta\models R^+$, since $(r_R, r_L\odot r_R)\in\theta$ for $(r_L, r_R)\in R$. By Proposition \ref{proposition:full_crossing_freest}, $\crossing_R M$ is freer or as free as $M/\theta$, through a surjective homomorphism $g$ with $g\circ\nu_{\crossing_R M} = \nu_{M/\theta}$. Hence $g\circ h_R\circ\nu_M = \pi_\theta\circ\nu_M$ and, as $\nu_M$ is surjective, $g\circ h_R = \pi_\theta$. Therefore $\ker h_R\subseteq\ker\pi_\theta = \theta$. In particular, if $\theta_R$ is $T$-aligned, then $\ker h_R\subseteq\theta_R\subseteq\ker h_R$, and conversely, if $\theta_R = \ker h_R$, then $\theta_R$ is $T$-aligned.

    \begin{enumerate}
        \item By Lemma \ref{lemma:principal_aligned}, $\Theta(b, a\odot b) = \theta_{\{(a,b)\}}$ is $T$-aligned, and $\crossing_{a\le b} M = \crossing_{\{(a,b)\}} M$ by Corollary \ref{corollary:freest_models} (2).

        \item Let $R = \{r_1,\dots,r_n\}$ with $r_i = (a_i, b_i)$, and let $\theta' = \theta_{\{r_1,\dots,r_{n-1}\}}$, which we may assume $T$-aligned by induction on $n$. Then $N = M/\theta'$ is a semilattice over $C$ with terms $T$ and $\pi_{\theta'}$ is interpretation-preserving, Lemma \ref{lemma:aligned}. Let $\theta''$ be the principal congruence on $N$ generated by the pair $(\pi_{\theta'}(b_n),\pi_{\theta'}(a_n\odot b_n))$, which is $T$-aligned by Lemma \ref{lemma:principal_aligned}. The kernel of $\pi_{\theta''}\circ\pi_{\theta'}$ is a congruence containing $\theta'$ and $(b_n, a_n\odot b_n)$, and, by the Correspondence Theorem, \cite{BurSan81}*{Chapter II, \S 6}, the preimage of $\theta''$ under $\pi_{\theta'}$ is contained in any such congruence, so $\ker(\pi_{\theta''}\circ\pi_{\theta'}) = \theta_R$. As $\pi_{\theta''}\circ\pi_{\theta'}$ is surjective and interpretation-preserving, $\theta_R$ is $T$-aligned by Lemma \ref{lemma:aligned}. Finally, since $x\,\theta\,y$ if and only if $x\,\theta\,x\odot y\,\theta\,y$, we have $\Theta(x,y) = \Theta(y, x\odot y)\vee\Theta(x, x\odot y)$, so every finitely generated congruence is $\theta_R$ for a finite $R$.

        \item By the previous identity, $\Theta(a,b) = \theta_R$ for $R = \{(b,a),(a,b)\}$, and $\crossing_R M = \crossing_{a\le b}\crossing_{b\le a} M$ by Definition \ref{definition:full_crossing_finite} and Corollary \ref{corollary:freest_models} (2). Apply (2).\qedhere
    \end{enumerate}
\end{proof}

Whether $\theta_R$ is $T$-aligned can be decided from the finite subsets of $R$.

\begin{corollary}\label{corollary:compactness}
    Let $M$ be a semilattice over $C$ with terms $T$ and let $R$ be a set of duples of $M$. Then $\theta_R$ is $T$-aligned, that is, $\crossing_R M = M/\theta_R$, if and only if, whenever $t\in T$ and $m\in M$ satisfy $(c\odot m, m)\in\theta_R$ for every $c\in C(t)$, there is a finite $R_0\subseteq R$ with $(c\odot m, m)\in\theta_{R_0}$ for every $c\in C(t)$.
\end{corollary}

\begin{proof}
    Recall that $\pi_\theta(x)\le\pi_\theta(y)$ if and only if $(x\odot y, y)\in\theta$. Assume the condition, let $t\in T$, and let $\pi_{\theta_R}(m)$ be an upper bound of the classes $\pi_{\theta_R}(c)$, $c\in C(t)$. Take $R_0$ as in the condition. By Theorem \ref{theorem:congruences}~(2), $\theta_{R_0}$ is $T$-aligned, so $(\nu_M(t)\odot m, m)\in\theta_{R_0}\subseteq\theta_R$. Conversely, let $\theta_R$ be $T$-aligned and let $t$ and $m$ be as in the condition. Then $(\nu_M(t)\odot m, m)\in\theta_R$. Since $\theta_R$ is the union of the directed family of the congruences $\theta_{R_0}$, $R_0\subseteq R$ finite, this pair lies in some $\theta_{R_0}$. For $c\in C(t)$, $c\odot m = (c\odot m)\odot m\ \theta_{R_0}\ (c\odot m)\odot\nu_M(t)\odot m = \nu_M(t)\odot m\ \theta_{R_0}\ m$.
\end{proof}

For an infinite set $R$ the inclusion $\theta_R\subseteq\ker h_R$ may be strict. Consider again the ascending chain $c_0<c_1<\cdots<c_\infty$ and the term $t$ with $C(t) = \{c_i\mid i\in\mathbb{N}\}$, and let $R = \{(c_i, c_0)\colon i\in\mathbb{N}\}$. Each $\Theta(c_0, c_i)$ identifies the elements of $[c_0, c_i]$, so $\theta_R$ has the two classes $\{c_i\mid i\in\mathbb{N}\}$ and $\{c_\infty\}$, that is, it is the congruence $\theta$ considered before Definition \ref{definition:aligned}. Thus $\theta_R$ is not $T$-aligned, its first class has no greatest element, and the condition of Corollary \ref{corollary:compactness} fails for $m=c_0$ and $t$, since a finite $R_0\subseteq R$ identifies $c_0$ with only finitely many $c_i$. On the other hand, a $T$-aligned congruence identifying $c_0$ with every $c_i$ also identifies $c_0$ with $\nu_M(t) = c_\infty$, so $\ker h_R = M\times M$ and $\crossing_R M$ is the one-element semilattice.

This result provides us with more context regarding full crossing. Indeed, the kernel of $\nu_M\colon F_{C;T}(\emptyset)\to M$ is a $T$-aligned congruence on $F_{C;T}(\emptyset)$, by Lemma \ref{lemma:aligned}, and $M = F_{C;T}(\emptyset)/\ker\nu_M$. For a finite set of duples, Theorem \ref{theorem:congruences} shows that full crossing is a succession of quotients by principal congruences, obtaining at each step a semilattice that is less free than the previous one, with the quotient projections as the successive homomorphisms. For an arbitrary set of duples $R$, full crossing computes the least $T$-aligned congruence containing $\theta_R$.

\section{Atom redundancy in infinite semilattices}\label{section:atom_redundancy}

One key element in the study of finite atomized semilattices is the concept of redundancy of atoms. Indeed, atomizations can be useful to have as they are set-theoretical based descriptions of semilattices that we can use to effectively perform several semilattice computations, but such descriptions are not practical if the number of atoms we require is too large. For example, in the case of the freest semilattice over a set of terms $T$, we have seen that the atom defined by every non-empty subset of $C$ is compatible with $F_{C;T}(\emptyset)$. This atomization has cardinality $2^{|C|}$, which can potentially be an infinite cardinal greater than $|C|$. However, we proved in Proposition \ref{proposition:freest_admits_atomization} that an atomization of $F_{C;T}(\emptyset)$ of cardinality $|C|$ always exists. Furthermore, we have already seen in Corollary \ref{corollary:cardinal_atomization_crossing} that performing the full crossing of a duple on a model with infinite atoms does not increase the cardinality of an atomization. This result immediately extends to crossing any finite set of duples.

In \cite{MarPol21}*{Definition 2.6}, the authors introduce the concept of redundant atoms, and show that redundant atoms can be added or removed from an atomization without changing the model. They then also show that there is a unique atomization of a model consisting exclusively of non-redundant atoms, and that every atomization for the model must include all of its non-redundant atoms.

Establishing analogous results in the infinite case is not straight-forward. In particular, if redundancy is introduced as it was done in the finite case, it is not always the case that removing redundant atoms does not change a model. As such, in this section we introduce two definitions of redundancy: redundancy and weak redundancy, Definition \ref{definition:redundancy}, and study their properties. These notions are equivalent in the finite case (and more generally in the complete case, as shown after Proposition \ref{proposition:redundancy_and_weak_redundancy}, and in semilattices with no infinite strictly ascending chain of terms, Corollary \ref{corollary:ascending_chain_least_atomization}), but not in general in the infinite one. Nonetheless, we are still able to show that, if $M$ is finitely generated, every atom of $M$ is a join of atoms non-redundant with $M$, see Theorem \ref{theorem:finitely_generated_join_redundant}.

We will now introduce some results exhibiting the differences between the finite and infinite case. Throughout this section, whenever we say that a semilattice $M$ over a set of constants $C$ is finitely generated, we are assuming that its set of terms $T$ is the set of finite subsets of $C$. In particular, and in accordance with Definition \ref{definition:duples}, the duples of $M$ will be the finitely generated duples over $C$. Similarly, whenever we say that $M$ is complete, we are assuming that its set of terms $T$ is the set of all non-empty subsets of $C$. A semilattice may be both finitely generated and complete, for instance any finite semilattice or the chain of Example \ref{example:infinite_chain} below; in that case, $T$ is the one fixed by the hypothesis named in each statement.

\begin{proposition}\label{proposition:atom_join_right_to_left}
    Let $M$ be an atomized semilattice atomized by a set of atoms $A$. If an atom $\phi$ is a union of atoms of $A$ different from $\phi$, both $A\setminus\{\phi\}$ and $A\cup\{\phi\}$ atomize $M$. More generally, if $B\subseteq\Omega(M)$ is a set of atoms such that every atom of $A$ is a union of atoms of $B$, then $B$ atomizes $M$.
\end{proposition}

\begin{proof}
    Let $r = (r_L, r_R)$ be a duple discriminated by $\phi$. Then, $\phi < r_L$ and $\phi\nless r_R$. Since $\phi<r_L$, there must be $c\in C(r_L)$ such that $\phi < c \leq r_L$. Given that $\phi$ is a union of the atoms in a set $S \subseteq A$ different from $\phi$, there must be at least one atom $\varphi \in S$ such that $\varphi < c \leq r_L$ and, because the upper constant segment of $\varphi$ is a subset of that of $\phi$, then $\varphi\nless r_R$. Since $\varphi < r_L$ and $\varphi\nless r_R$, $\varphi$ discriminates $r$. It follows that for any duple discriminated by $\phi$ there is an atom of $S$ that discriminates the duple. Hence every duple discriminated by $\phi$ is negative in $M$, so $\phi\in\Omega(M)$ and $A\cup\{\phi\}$ atomizes $M$, as shown after Definition \ref{definition:compatible}; and every negative duple discriminated by $\phi$ is still discriminated by an atom of $A\setminus\{\phi\}$, so $A\setminus\{\phi\}$ satisfies \eqref{axiom:AS3}. Finally, \eqref{axiom:AS6} holds for $A\setminus\{\phi\}$, since every constant in $U^c(\phi)$ is in the upper constant segment of some atom of $S$.

    For the second statement, the atoms of $B$ are compatible with $M$, so they discriminate no duple positive in $M$. If $M\models r^-$, some $\phi\in A$ discriminates $r$, and, as above, an atom $\varphi\in B$ with $U^c(\varphi)\subseteq U^c(\phi)$ discriminates $r$. Thus $B$ satisfies \eqref{axiom:AS3}. It also satisfies \eqref{axiom:AS6}, since the upper constant segments of the atoms of $B$ cover those of the atoms of $A$, which cover $C$. Hence $B$ atomizes $M$, Proposition \ref{proposition:spawned}.
\end{proof}

In the finite case, Proposition \ref{proposition:atom_join_right_to_left} is an equivalence. We will show that this is still the case in the complete case. For that, we need to introduce pinning terms, following the definition from \cite{MarPol21}.

\begin{definition}
    Let $M$ be a semilattice over a set of constants $C$ with terms $T$. Let $\phi$ be an atom, $\phi\ne\ominus_C$. We say that the term $T_\phi = \odot \{c : c\in C\setminus U^c(\phi) \}$ is the \emph{pinning term of $\phi$}. We say that the pinning term exists for $M$ if $T_\phi$ is in $T$. If the pinning term exists, for each $c\in U^c(\phi)$, we say that $(c,T_\phi)$ is a \emph{pinning duple} of $\phi$.
\end{definition}

More generally, we say that a term $s\in T$ \emph{pins} $\phi$ if $C(s)\cap U^c(\phi) = \emptyset$ and $M\models c\le s$ for every $c\in C\setminus U^c(\phi)$. If $T_\phi\in T$, it pins $\phi$, and every term pinning $\phi$ has the same value as $T_\phi$.

One key difference between the finite and infinite setting is that in the latter pinning terms do not necessarily exist. For example, in a finitely generated atomized semilattice over an infinite set of constants, the pinning term of an atom with a finite upper constant segment, as defined in \cite{MarPol21}, does not exist, although, with the more general definition given above, other terms may still pin the atom. For instance, if $M$ is finitely generated over $C = \{c_0, c_1, c_2,\dots\}$ with $c_0 > c_1 > c_2 > \cdots$ and $\phi_i :< \{c_0,\dots,c_i\}$, then $T_{\phi_i}$ is not a term, but $c_{i+1}$ pins $\phi_i$. However, they always exist in complete semilattices, allowing us to extend Proposition \ref{proposition:atom_join_right_to_left} to an equivalence.

\begin{theorem}\label{theorem:complete_semilattice_atom_join}
    Let $M$ be a complete atomized semilattice. Assume that $M$ is atomized by a set of atoms $A$. Then, the sets $A\setminus\{\phi\}$ and $A\cup\{\phi\}$ both atomize $M$ if and only if $\phi$ is a union of atoms of $A$ different from $\phi$.
\end{theorem}

\begin{proof}
    Right to left follows from Proposition \ref{proposition:atom_join_right_to_left}. To prove left to right, assume that $A\setminus\{\phi\}$ and $A\cup\{\phi\}$ both atomize $M$. If $\phi=\ominus_C$, \eqref{axiom:AS6} for $A\setminus\{\phi\}$ and \eqref{axiom:AS4} give that every constant is in the upper constant segment of some atom of $A\setminus\{\ominus_C\}$, so $\ominus_C$ is the union of the atoms of $A\setminus\{\ominus_C\}$. Assume now that $\phi\ne\ominus_C$. Since $M$ is complete, the pinning term $T_\phi$ exists, and by definition $\phi\nless T_\phi$. Then, $\phi$ discriminates the pinning duple $(c, T_\phi)$. In fact, it does so for every $c\in U^c(\phi)$.

    Since $A\setminus\{\phi\}$ and $A\cup\{\phi\}$ both atomize $M$, every pinning duple $(c, T_\phi)$ is negative in $M$, with or without $\phi$, so every such pinning duple must be discriminated by some other atom. Hence, for every $c\in U^c(\phi)$, there is an atom $\psi_c$ in $A\setminus\{\phi\}$ discriminating $(c, T_\phi)$, and then $\phi = \join_{c\in U^c(\phi)} \psi_c$.
\end{proof}

Then, if the semilattice $M$ is complete, for each duple $r$ discriminated by an atom $\phi$ that can be added or removed from a model, there is another atom $\psi_r$ that discriminates $r$ and satisfies
$\psi_r \join \phi = \phi$. This, in particular, holds for semilattices over a finite set of constants, as they are always complete.

In the finitely generated case, since pinning terms may not exist, we can instead use a similar argument with any term generated by a finite set of constants that does not intersect the upper constant segment of the atom, thus proving that if an atom $\phi$ that can be added or removed to a model discriminates a duple $r$, there is another atom $\psi_r$ in the model discriminating such duple. However, it may not be the case that $\phi$ is wider than $\psi_r$. Consequently, we cannot guarantee that $\phi$ can be obtained as a join of atoms. Indeed, consider the following example.

\begin{example}\label{example:infinite_chain}
    Suppose that we have a model over $C = \{c_0, c_1,\dots,c_\infty\}$ where $c_i < c_{i+1}$ and $c_i < c_\infty$, for all $i\ge 0$. Since $c_{i+1}\not\le c_i$, we need an atom discriminating this duple: an atom in $c_{i+1}$ not in $c_i$. This atom should also contain $c_j$, for every $j\ge i+1$, as well as $c_\infty$. We could then think that this model is atomized by the set of atoms $A = \{\phi_0, \phi_1,\dots\}$ where $\phi_i :< \{c_\infty\}\cup\{c_j\mid j\ge i\}$.

    This atomizes such a model, with a caveat. We have not specified the set of admissible terms in our model. If the model is finitely generated, there is no ambiguity. However, if the model is complete, this atomization is not able to discriminate between $c_\infty$ and any infinite join of the remaining constants. Therefore, adding the atom $\phi_\infty$ whose upper constant segment is $U^c(\phi_\infty) = \{c_\infty\}$ would change the model if the semilattice is complete, but not if it is finitely generated. Notice further that this atom $\phi_\infty$ cannot be written as a join of the atoms in our atomization.
\end{example}

We now show that an atom being a join of atoms in a complete atomized semilattice is independent of the chosen atomization.

\begin{proposition}\label{proposition:redundancy_if_complete}
    Let $M$ be a complete semilattice over a set of constants $C$ atomized by a set $A$ of atoms. Let $\phi$ be an atom over $C$. Then $\phi$ is a union of atoms of $A$ different from $\phi$ if and only if $\phi$ is a union of atoms of $\Omega(M)$ different from $\phi$.
\end{proposition}

\begin{proof}
    Left to right is straightforward. Right to left, assume that $\phi = \join_i \psi_i$, with $\phi\ne\psi_i\in\Omega(M)$. If $\psi_i$ is not an atom in $A$, and given that $\psi_i\in\Omega(M)$, $A\cup\{\psi_i\}$ still atomizes $M$. Furthermore, $A\setminus\{\psi_i\} = A$ also atomizes $M$. Then by Theorem \ref{theorem:complete_semilattice_atom_join}, $\psi_i = \join_j\varphi_j$, where each $\varphi_j\in A$, and $\varphi_j\ne\phi$ since $U^c(\varphi_j)\subseteq U^c(\psi_i)\subsetneq U^c(\phi)$. By doing this with every $\psi_i$ not in $A$, we get to write $\phi$ as a union of atoms in $A$.
\end{proof}

Notice that this may fail to hold in non-complete semilattices. For instance, let $M = F_C(\emptyset)$ be the finitely generated model atomized by the set $A$ of atoms $\phi_{ij}$ of Example \ref{example:every_atom_redundant} below, and let $\phi:<\{c_0, c_1\}$. Every atom over $C$ is compatible with $M$, Proposition \ref{proposition:all_atoms_fit_freest}, so $\phi = \psi_{c_0}\join\psi_{c_1}$ is a union of atoms of $\Omega(M)$ different from $\phi$, where $\psi_c:<\{c\}$. However, $\phi$ is not a union of atoms of $A$, since every atom of $A$ has an infinite upper constant segment.

So now we raise the question of what properties an atom needs to fulfil in order for it to be removable from a finitely generated semilattice without changing the model. We have the following.

\begin{theorem}\label{theorem:weakly_redundant_removable}
    Let $M$ be an atomized semilattice over a set of constants $C$ with terms $T$. Let $A$ be a set of atoms atomizing $M$ and let $\phi$ be an atom. The sets $A\setminus\{\phi\}$ and $A\cup\{\phi\}$ both atomize $M$ if and only if for every constant $c\in U^c(\phi)$ and every $t\in T$ with $C(t)\subseteq C\setminus U^c(\phi)$, there is an atom $\psi\in A$ such that $\psi\ne\phi$, $\psi < c$, $\psi\nless t$, and for every constant $c\in U^c(\phi)$ there is an atom $\psi\in A$ such that $\psi\ne\phi$, $\psi < c$.
\end{theorem}

\begin{proof}
    Every duple discriminated by $\phi$ is of form $(s,t)$ where $\phi < s$ and $C(t)\subseteq C\setminus U^c(\phi)$. Since $\phi < s$ implies that there exists a constant $c\in C(s)$ such that $\phi < c$, $\phi$ also discriminates the duple $(c,t)$. Reciprocally, an atom discriminating $(c,t)$ with $c\in C(s)$ must also discriminate $(s,t)$. We thus only need to consider duples of the form $(c,t)$, with $c\in U^c(\phi)$ and $C(t)\subseteq C\setminus U^c(\phi)$.

    Then, $A\setminus\{\phi\}$ and $A\cup\{\phi\}$ both satisfy \eqref{axiom:AS3} for $M$ if and only if every such duple $(c,t)$ is also discriminated by another atom in $A$; that is, an atom $\psi\in A$ such that $\psi\ne\phi$, $\psi < c$, $\psi\nless t$.

    Axioms \eqref{axiom:AS1}, \eqref{axiom:AS2}, \eqref{axiom:AS4} and \eqref{axiom:AS5} hold for any subset of $\Omega(M)$, see the proof of Theorem \ref{theorem:every_semilattice_atomizable}, and \eqref{axiom:AS6} holds for $A\cup\{\phi\}$ because it holds for $A$. By \eqref{axiom:AS4}, an atom is lower than an element of $M$ if and only if it is lower than one of its constants, so \eqref{axiom:AS6} holds for $A\setminus\{\phi\}$ if and only if every constant is greater than some atom of $A\setminus\{\phi\}$. A constant not in $U^c(\phi)$ is greater than some atom of $A$, which cannot be $\phi$. Hence \eqref{axiom:AS6} holds for $A\setminus\{\phi\}$ if and only if for every $c\in U^c(\phi)$ there is an atom $\psi\in A$ such that $\psi\ne\phi$, $\psi<c$.
\end{proof}

The results above motivate the introduction of the following definition.

\begin{definition}\label{definition:redundancy}
    Let $M$ be an atomized semilattice over a set of constants $C$ with terms $T$ and let $\phi$ be an atom. We say that $\phi$ is:
    \begin{itemize}
        \item \emph{redundant with $M$} if $\phi$ is a union of atoms of $\Omega(M)$ different from $\phi$;
        \item \emph{weakly redundant with $M$} if for every constant $c \in C$ and every term $t\in T$ such that $\phi$ discriminates $(c, t)$ there is an atom $\psi\in\Omega(M)$ different from $\phi$ that also discriminates $(c, t)$, and for every constant $c\in U^c(\phi)$ there is an atom $\psi\in\Omega(M)$ different from $\phi$ with $\psi < c$.
        \item \emph{non-redundant with $M$} if $\phi\in\Omega(M)$ and $\phi$ is not redundant with $M$.
        \item \emph{non-weakly-redundant with $M$} if $\phi\in\Omega(M)$ and $\phi$ is not weakly redundant with $M$.
    \end{itemize}
    The sets of atoms redundant with $M$, weakly redundant with $M$, non-redundant with $M$, and non-weakly-redundant with $M$ are respectively denoted by $\R(M)$, $\WR(M)$, $\NR(M)$, and $\NWR(M)$.
\end{definition}

The second condition of weak redundancy only matters for the zero atom. If $\phi\ne\ominus_C$, let $c\in U^c(\phi)$ and let $d$ be a constant not in $U^c(\phi)$; then $\phi$ discriminates $(c,d)$, and an atom different from $\phi$ discriminating $(c,d)$ is lower than $c$. The zero atom $\ominus_C$ discriminates no duple, so it is weakly redundant with $M$ if and only if the upper constant segments of the atoms of $\Omega(M)$ other than $\ominus_C$ cover $C$, that is, if and only if $\ominus_C$ is redundant with $M$. As in the finite case, \cite{MarPol21}*{Theorem 6}, this happens if and only if $M$ has no least element. A least element of $M$, if it exists, is the value of a constant: it is $\nu_M(t)$ for some $t\in T$, and every $c\in C(t)$ satisfies $\nu_M(t)\le c\le\nu_M(t)$. So it suffices to look at constants. If a constant $c$ is not a least element of $M$, there is a duple $(c,t)$ negative in $M$, and an atom of $\Omega(M)$ discriminating it is lower than $c$ and different from $\ominus_C$. If $c$ is a least element of $M$, an atom $\phi\ne\ominus_C$ with $\phi<c$ discriminates the positive duple $(c,d)$ for any constant $d\notin U^c(\phi)$, so it is not compatible with $M$.

\begin{lemma}\label{lemma:redundant_implies_in_model}
    Let $M$ be an atomized semilattice over a set of constants $C$ and let $\phi$ be an atom. If $\phi$ is redundant or weakly redundant with $M$, $\phi\in\Omega(M)$.
\end{lemma}

\begin{proof}
    Assume first that $\phi$ is redundant. Then $\phi$ is a union of atoms of $\Omega(M)$ and thus, by Corollary \ref{corollary:join_still_in_model}, it is an element of $\Omega(M)$.

    Assume now that $\phi$ is weakly redundant, and let $r=(r_L, r_R)$ be a duple discriminated by $\phi$. Then $\phi$ also discriminates $(c, r_R)$ for some $c\in C(r_L)$, and by weak redundancy an atom $\psi\in\Omega(M)$ discriminates $(c, r_R)$, hence also $r$. Since $\psi$ is compatible with $M$, $M\models r^-$. Therefore, $\phi\in\Omega(M)$.
\end{proof}

Notice that, if $M$ is complete and atomized by $A$, an atom $\phi$ is redundant if and only if it is equal to a union of atoms of $A$ different from $\phi$, as a consequence of Proposition \ref{proposition:redundancy_if_complete}. Notice further that it may be the case that a weakly redundant atom cannot be removed from an atomization $A$ of a model $M$ without changing the model because, although there are other atoms in $\Omega(M)$ that discriminate every duple discriminated by $\phi$, such atoms are not necessarily in $A$. On the other hand, a non-weakly-redundant atom can never be removed from any atomization, as shown in Corollary \ref{corollary:every_atomization_contains_non-weekly_redundant} below.

\begin{proposition}
    Let $M$ be an atomized semilattice over a set of constants $C$ with terms $T$. An atom $\phi\in\Omega(M)$ is non-weakly-redundant if and only if there exists a duple $r$ of $M$ such that $M\models r^-$ and such that $\phi$ is the only atom compatible with $M$ discriminating $r$, or there exists a constant $c$ such that $\phi$ is the only atom compatible with $M$ lower than $c$.
\end{proposition}

\begin{proof}
    Assume first that $\phi\in\NWR(M)$. By definition, one of the two conditions of weak redundancy fails. If the first one fails, there exists $c\in U^c(\phi)$ and $t\in T$ with $C(t)\subseteq C\setminus U^c(\phi)$ such that there is no atom different from $\phi$ discriminating $r = (c,t)$. However, $\phi\in\Omega(M)$ discriminates $r$, thus $M\models r^-$. If the second one fails, there is a constant $c\in U^c(\phi)$ such that $\phi$ is the only atom of $\Omega(M)$ lower than $c$.

    Assume now that there is a duple $r$ discriminated by an atom $\phi$, and only by $\phi$, or a constant $c$ such that $\phi$ is the only atom of $\Omega(M)$ lower than $c$. Then $\Omega(M)\setminus\{\phi\}$ does not atomize $M$, as it fails \eqref{axiom:AS3} or \eqref{axiom:AS6} respectively, and by Theorem \ref{theorem:weakly_redundant_removable}, $\phi\in\NWR(M)$.
\end{proof}

For $\phi\ne\ominus_C$ the second alternative implies the first, as shown after Definition \ref{definition:redundancy}. For $\ominus_C$, which discriminates no duple, only the second alternative can occur, and it says that $c$ is a least element of $M$.

Consequently, we have the following.

\begin{corollary}\label{corollary:every_atomization_contains_non-weekly_redundant}
    Let $C$ be a set of constants and let $M$ be a semilattice over $C$. Then $\NWR(M)$ is the intersection of all atomizations of $M$; in particular, $\NWR(M)\subseteq A$ for every atomization $A$ of $M$. Consequently, $M$ has a least atomization if and only if $\NWR(M)$ atomizes $M$, and in that case the least atomization is $\NWR(M)$.
\end{corollary}

\begin{proof}
    Let $A$ atomize $M$ and $\phi\in\NWR(M)$. If $\phi\notin A$, then $A\subseteq\Omega(M)\setminus\{\phi\}$, so $\Omega(M)\setminus\{\phi\}$ atomizes $M$, as shown after Definition \ref{definition:compatible}, against the previous proposition. Conversely, let $\phi$ belong to every atomization of $M$. Then $\phi\in\Omega(M)$, since $\Omega(M)$ atomizes $M$. If $\phi$ were weakly redundant, the conditions of Theorem \ref{theorem:weakly_redundant_removable} would hold for $A = \Omega(M)$, and $\Omega(M)\setminus\{\phi\}$ would be an atomization of $M$ not containing $\phi$. Hence $\phi\in\NWR(M)$. The last statement follows, since an atomization contained in every atomization is contained in their intersection.
\end{proof}

We have the following relationship between redundancy and weak redundancy.

\begin{proposition}\label{proposition:redundancy_and_weak_redundancy}
    Let $C$ be a constant set, and let $M$ be a semilattice over $C$ with terms $T$. Let $\phi\in\Omega(M)$ be an atom.
    \begin{enumerate}
        \item If $\phi$ is redundant, it is weakly redundant.
        \item If $\phi$ is weakly redundant and some term of $T$ pins $\phi$, it is redundant.
    \end{enumerate}
\end{proposition}

\begin{proof}
    Assume first that $\phi$ is redundant. Then $\phi = \join_i \psi_i$, where $\psi_i$ are atoms in $\Omega(M)$ different from $\phi$. In particular, for every $c\in U^c(\phi)$ there is some $i$ such that $\psi_i < c$, which is the second condition of weak redundancy. Then, if we have a duple $(c,t)$ with $c\in U^c(\phi)$ and $C(t)\subseteq C\setminus U^c(\phi)$, there is some $i$ such that $\psi_i < c$. Furthermore, $U^c(\psi_i)\subsetneq U^c(\phi)$, thus $\psi_i\ne \phi$ and $\psi_i\nless t$.

    Assume now that $\phi$ is weakly redundant and that $s\in T$ pins $\phi$. Then $C(s)\subseteq C\setminus U^c(\phi)$, so for every $c\in U^c(\phi)$ there is an atom $\psi_c$ such that $\psi_c\ne \phi$, $\psi_c < c$, $\psi_c\nless s$. This last inequality implies that $U^c(\psi_c)\subseteq U^c(\phi)$: if $d\in U^c(\psi_c)\setminus U^c(\phi)$, then $\psi_c < d\le s$. We then have that $\phi=\join_{c\in U^c(\phi)}\psi_c$, exhibiting that $\phi$ is redundant.
\end{proof}

In particular, if $M$ is a complete semilattice (and particularly if $C$ is finite) pinning terms exist for every atom other than $\ominus_C$, so an atom is redundant in $M$ if and only if it is weakly redundant; for $\ominus_C$ this holds in every semilattice, as shown after Definition \ref{definition:redundancy}. Thus, if $M$ is complete, $\NR(M) = \NWR(M)$, so, by Corollary \ref{corollary:every_atomization_contains_non-weekly_redundant}, $\NR(M)$ is contained in every atomization of $M$, and $M$ has a least atomization if and only if $\NR(M)$ atomizes $M$, as in the finite case.

Pinning also tells exactly which non-redundant atoms are weakly redundant.

\begin{proposition}\label{proposition:NWR_pinned}
    Let $M$ be a semilattice over a set of constants $C$ with terms $T$, and let $\phi\in\Omega(M)$, $\phi\ne\ominus_C$. Then $\phi$ is non-weakly-redundant if and only if $\phi$ is non-redundant and some term of $T$ pins $\phi$.
\end{proposition}

\begin{proof}
    If $\phi$ is non-redundant and some term of $T$ pins $\phi$, then $\phi$ is not weakly redundant, by Proposition \ref{proposition:redundancy_and_weak_redundancy}~(2), so $\phi\in\NWR(M)$. Conversely, let $\phi\in\NWR(M)$. Then $\phi\in\NR(M)$, by Proposition \ref{proposition:redundancy_and_weak_redundancy}~(1). Since $\phi\ne\ominus_C$, the first condition of weak redundancy fails, as shown after Definition \ref{definition:redundancy}, so there are $c\in C$ and $t\in T$ such that $\phi$ discriminates $(c,t)$ and no other atom of $\Omega(M)$ does. Then $C(t)\cap U^c(\phi)=\emptyset$. Let $K=\{d\in C\colon M\models (d\not\le t)\}$. Since $\phi\in\Omega(M)$, the duple $(c,t)$ is negative in $M$, so $c\in K$, while no constant of $t$ is in $K$. As in the proof of Lemma \ref{lemma:compatible_atom_exists}, the atom $\phi' :< K$ is compatible with $M$ and discriminates $(c,t)$. Hence $\phi'=\phi$, that is, $U^c(\phi)=K$, and $M\models d\le t$ for every $d\in C\setminus U^c(\phi)$. Thus $t$ pins $\phi$.
\end{proof}

Hence $\NR(M)\setminus\NWR(M)$ consists of the non-redundant atoms other than $\ominus_C$ that no term of $T$ pins.

Without a pinning term, the converse of Proposition \ref{proposition:redundancy_and_weak_redundancy} (1) may fail. In the finitely generated reading of Example \ref{example:infinite_chain}, $\phi_\infty$ is weakly redundant, since a duple $(c_\infty, t)$ discriminated by $\phi_\infty$ has $t$ finite and is also discriminated by $\phi_{k+1}$, where $k$ is the largest index of a constant in $t$, and $\phi_0 < c_\infty$. However, $\phi_\infty$ is not redundant, since it is the only atom whose upper constant segment is contained in $\{c_\infty\}$. Moreover, $\Omega(M) = A\cup\{\phi_\infty\}$, since a compatible atom lower than $c_i$ is lower than every constant above $c_i$, and every atom of $\Omega(M)$ is non-redundant, since a compatible atom whose upper constant segment is strictly contained in $U^c(\phi_i)$ is not lower than $c_i$. On the other hand, for $i\ge 1$, $\phi_i$ is the only atom of $\Omega(M)$ discriminating $(c_i, c_{i-1})$, and $\phi_0 = \ominus_C$ is non-weakly-redundant because $c_0$ is the least element of $M$, as shown after Definition \ref{definition:redundancy}. Hence $A = \NWR(M)$ is the least atomization of $M$, Corollary \ref{corollary:every_atomization_contains_non-weekly_redundant}, and it omits the non-redundant atom $\phi_\infty$.

However, if $M$ is not complete, even when $\NR(M)$ atomizes $M$ there may be atomizations that do not contain any non-redundant atoms, thus $\NWR(M) = \emptyset$; and there may be cases where $M$ is complete and has no non-redundant atoms, as the following example shows.

\begin{example}\label{example:every_atom_redundant}
    Let $C = \{c_0, c_1, c_2,\dots\}$. Consider the atomization provided by atoms $\phi_{ij} :< \{c_k \mid k = j \mod 2^i\}$, where $i = 0,1,2,\dots$ and $j\in \{0,1,\dots, 2^i - 1\}$. 

    If we consider $M$ a finitely generated model over $C$ with this atomization, it is the freest model. Indeed, assume that $r = (r_L, r_R)$ is a finitely generated duple such that $C(r_L)\not\subseteq C(r_R)$. Let $j\ge 0$ be such that $c_j\in C(r_L)\setminus C(r_R)$ and let $i$ be such that $j < 2^i$ and $\max\{k\mid c_k\in C(r_R)\} < 2^i$. Then, $\phi_{i j}$ discriminates $r$, as $c_j\in U^c(\phi_{ij})$ and every other constant in the upper constant segment has an index greater than that of any constant in $C(r_R)$.

    Notice, however, that this atomization consists exclusively of redundant (and thus weakly redundant) atoms, as $\phi_{ij} = \phi_{i + 1, j}\join \phi_{i + 1, j + 2^i}$. Since every atomization must contain every non-weakly-redundant atom (Corollary \ref{corollary:every_atomization_contains_non-weekly_redundant}), $\NWR\big(F_C(\emptyset)\big) = \emptyset$, so, by the same corollary, $F_C(\emptyset)$ has no least atomization. Nonetheless, $F_C(\emptyset)$ is atomized by $\NR\big(F_C(\emptyset)\big)$, as the atomization provided in Proposition \ref{proposition:freest_admits_atomization} consists exclusively of non-redundant atoms.

    On the other hand, we can also consider the complete model $M$ over $C$ spawned by the same atoms. This model is not the freest complete semilattice over $C$, as the duple $(c_0,\odot_{i > 0} c_i)$ is clearly positive. Due to a similar reasoning, it is still the case that $\NWR(M) = \emptyset$, but given that in complete models an atom is non-weakly-redundant if and only if it is non-redundant, $\NR(M)=\emptyset$. This example thus provides a pathological example of a complete model in which every atom is redundant.
\end{example}

Nonetheless, one thing that we can show is that crossing a finite number of duples starting from the freest model does produce models that can be atomized by non-(weakly)-redundant atoms. To do so, let us first introduce the following terminology.

\begin{definition}
    Let $C$ be a set of constants and let $\phi$ be an atom over $C$. We say that $\phi$ is:
    \begin{enumerate}
        \item \emph{finite} if $U^c(\phi)$ is finite;
        \item \emph{infinite} if $U^c(\phi)$ is infinite;
        \item \emph{cofinite} if $C\setminus U^c(\phi)$ is finite;
        \item \emph{coinfinite} if $C\setminus U^c(\phi)$ is infinite.
    \end{enumerate}
\end{definition}

We then start with the following.

\begin{proposition}\label{proposition:finite_constant_segments_non_redundant}
    Let $M$ be a semilattice. Assume that $M$ can be atomized by a set of finite atoms $A$. Then $M$ is atomized by $\NR(M)$.
\end{proposition}

\begin{proof}
    Let $\phi\in A$ be an atom. If $\phi$ is redundant, we can write it as a join of atoms different from $\phi$, $\phi=\join_i \psi_i$, $\psi_i\in\Omega(M)$. If all the $\psi_i$ are non-redundant, $\phi$ can be written as a join of non-redundant atoms. Otherwise, we can perform the same step with the redundant $\psi_i$. Since all the involved atoms are finite, this process must finish in a finite number of steps, and we thus deduce that $\phi$ can be written as a join of atoms non-redundant with $M$.

    We have thus shown that every atom in an atomization of $M$ is in $\NR(M)$ or can be written as a join of atoms in $\NR(M)$, which in particular implies that $\NR(M)$ is an atomization for $M$, Proposition \ref{proposition:atom_join_right_to_left}.
\end{proof}

\begin{corollary}\label{corollary:finite-implies-finite-after-crossing}
    Let $M$ be a semilattice atomized by a set of finite atoms $A$. Let $r$ be a duple of $M$. Then $\crossing_r M$ can also be atomized by a set of finite atoms. In particular, it is atomized by $\NR(\crossing_r M)$.
\end{corollary}

\begin{proof}
    The full crossing $\crossing_r M$ is atomized by $\big(A\setminus\dis_M(r)\big)\cup \big(\dis_M(r)\join L_M^a(r_R)\big)$. Since the atoms of $A$ are finite, so are the atoms in $A\setminus\dis_M(r)$. On the other hand, the atoms in $\dis_M(r)\join L_M^a(r_R)$ are the join of two atoms in $A$, so they again must be finite. We deduce that $\crossing_r M$ is atomized by an atomization consisting of finite atoms. It then follows by Proposition \ref{proposition:finite_constant_segments_non_redundant} that $\crossing_r M$ is atomized by $\NR(\crossing_r M)$.
\end{proof}

We can prove the following result, telling us that models of finite theories can be atomized by their non-redundant atoms. By finite theory we mean a theory consisting of a finite set of duples naming any terms, where terms can be finite or infinite.

\begin{corollary}\label{corollary:finite_theory_non_redundant}
    Let $C$ be a set of constants, $T$ a set of terms over $C$ and $R$ a finite set of duples of $T$. Then $F_{C;T}(R)$ is atomized by its non-redundant atoms.
\end{corollary}

\begin{proof}
    By Proposition \ref{proposition:freest_admits_atomization}, $F_{C;T}(\emptyset)$ admits an atomization consisting exclusively of finite atoms. Write $R=\{r_1,\dots,r_n\}$. Every semilattice with terms $T$ models $\Th_0^+(F_{C;T}(\emptyset))$, Proposition \ref{proposition:atoms_and_freest_model}~(1), so, by Theorem \ref{theorem:full_crossing_freest} applied as in the proof of Corollary \ref{corollary:commutativity}, $\crossing_{r_n}\cdots\crossing_{r_1}F_{C;T}(\emptyset) = F_{C;T}(R)$. By Corollary \ref{corollary:finite-implies-finite-after-crossing}, applied $n$ times, $F_{C;T}(R)$ is atomized by a set of finite atoms, and hence by its non-redundant atoms, Proposition \ref{proposition:finite_constant_segments_non_redundant}.
\end{proof}

Notice that, in the case of complete models, this equivalently tells us that the model is atomized by its non-weakly-redundant atoms, which then form its least atomization, Corollary \ref{corollary:every_atomization_contains_non-weekly_redundant}.

\begin{figure}
    \begin{tikzpicture}
        \draw (-1.5, 0) ellipse (2.2 cm and 1.5 cm);
        \draw (1.5, 0) ellipse (2.2 cm and 1.5 cm);
        \node (redundant) at (-2.2, 0) {$\R(M)$};
        \node [fill=white] (weakly-redundant) at (-1.5, -1.5) {$\WR(M)$};
        \node [fill=white] (non-redundant) at (1.5, -1.5) {$\NR(M)$};
        \node (non-weakly-redundant) at (2.4, 0) {$\NWR(M)$};
    \end{tikzpicture}
    \caption{\label{figure:decomposition_general}A representation of the decomposition of the set of atoms $\Omega(M)$ compatible with a semilattice $M$, following Lemma \ref{lemma:redundant_implies_in_model}, Proposition \ref{proposition:redundancy_and_weak_redundancy} and Proposition \ref{proposition:compatible_atoms_decomposition}.}
\end{figure}
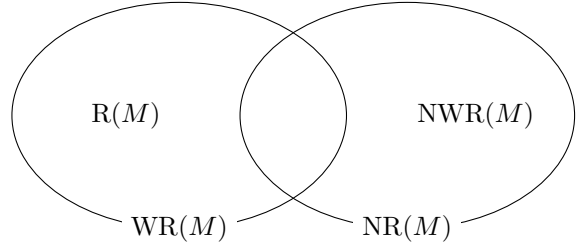

\begin{proposition}\label{proposition:compatible_atoms_decomposition}
    Let $M$ be an atomized semilattice over a set of constants $C$.
    \begin{enumerate}
        \item $\Omega(M)$ is the union of the set of redundant atoms with $M$ and the set of atoms non-redundant with $M.$
        \item $\Omega(M)$ is the union of the set of weakly redundant atoms with $M$ and the set of atoms non-redundant with $M,$ see Figure \ref{figure:decomposition_general}.
    \end{enumerate}
\end{proposition}

\begin{proof}
    (1) follows from the definition of non-redundant atoms. (2) follows from (1), Lemma \ref{lemma:redundant_implies_in_model} and Proposition \ref{proposition:redundancy_and_weak_redundancy} (1).
\end{proof}

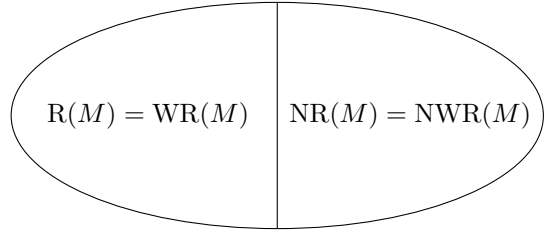
\begin{figure}
    \begin{tikzpicture}
        \draw (0, 0) ellipse (3.5 cm and 1.5 cm);
        \draw (0, -1.5) -- (0, 1.5);
        \node (redundant) at (-1.7, 0) {$\R(M) = \WR(M)$};
        \node (non-redundant) at (1.75, 0) {$\NR(M) = \NWR(M)$};
    \end{tikzpicture}
    \caption{A representation of the decomposition of the set of atoms $\Omega(M)$ compatible with a complete semilattice $M$, or with a semilattice whose terms admit no infinite ascending chain distinct in $M$, following Lemma \ref{lemma:redundant_implies_in_model}, Proposition \ref{proposition:redundancy_and_weak_redundancy}, Corollary \ref{corollary:ascending_chain_least_atomization} and Proposition \ref{proposition:compatible_atoms_decomposition}. In particular, this includes semilattices over a finite set of constants.}
\end{figure}

We now introduce some results regarding the completion of semilattices, Definition \ref{definition:completion}. 

\begin{proposition}\label{proposition:completion_atomization_base}
    Let $S$ be a semilattice over a constant set $C$ with terms $T$, and let $M$ be a complete semilattice over $C$.
    \begin{enumerate}
        \item $M$ is a completion of $S$ if and only if $S$ and $M$ have the same positive duples over $T$.
        \item If $M$ is a completion of $S$, any atomization of $M$ is also an atomization of $S$.
    \end{enumerate}
\end{proposition}

\begin{proof}
    (1) Assume first that $M$ is a completion of $S$, and let $i\colon S\to M$ be the injective, interpretation-preserving homomorphism. Let $r = (r_L, r_R)$ be a duple of $S$. Then $S\models r^+$ if and only if $\nu_S(r_R) = \nu_S(r_L)\odot \nu_S(r_R)$. Since $i$ is an injective homomorphism with $i(\nu_S(t)) = \nu_M(t)$ for $t\in T$, this holds if and only if $\nu_M(r_R) = \nu_M(r_L)\odot \nu_M(r_R)$, that is, if and only if $M\models r^+$. Conversely, assume that $S$ and $M$ have the same positive duples over $T$. Every element of $S$ is $\nu_S(t)$ for some $t\in T$, Definition \ref{definition:admissible_terms}. Define $i\colon S\to M$ by $i(\nu_S(t)) = \nu_M(t)$ for $t\in T$. For $t,s\in T$, $\nu_S(t) = \nu_S(s)$ if and only if $S\models (t\le s)$ and $S\models (s\le t)$, which holds if and only if the same holds in $M$, that is, if and only if $\nu_M(t) = \nu_M(s)$. So $i$ is well defined and injective. It is a homomorphism because $t\cup s\in T$ for $t,s\in T$, so that $i(\nu_S(t)\odot\nu_S(s)) = i(\nu_S(t\cup s)) = \nu_M(t\cup s) = \nu_M(t)\odot\nu_M(s)$. By construction, $i$ is interpretation-preserving.

    (2) By (1), the duples of $S$ are positive if and only if they are also positive in $M$, so any atomization of $M$ is also an atomization of $S$.
\end{proof}

However, different atomizations for $S$ may produce different complete models. We can then introduce the following.

\begin{definition}
    Let $S$ be a semilattice over a set of constants $C$ with terms $T$. We say that two atomizations $A_1$ and $A_2$ of $S$ are \emph{completely equivalent} if the complete semilattices they spawn have the same positive theories over the set of every duple over $C$, with finite or infinite terms.
\end{definition}

Thus, from this point of view, a semilattice $S$ has as many classes of completely equivalent atomizations as different completions. Furthermore, completely equivalent atomizations of a semilattice can still be quite different. Nonetheless, we can use non-redundant atoms to establish some commonalities.

\begin{proposition}
    Let $S$ be a semilattice and let $A_1$ and $A_2$ be completely equivalent atomizations of $S$. Then, $A_1\cap \NR(S) = A_2 \cap \NR(S)$ and $A_1\cap \NWR(S) = A_2\cap\NWR(S)$.
\end{proposition}

\begin{proof}
    Let $M$ be the complete model atomized by $A_1$, a completion of $S$ by Proposition \ref{proposition:completion_atomization_base}~(1), since $A_1$ also atomizes $S$. Since $A_1$ and $A_2$ are completely equivalent atomizations, $A_2$ also atomizes $M$.

    Let $\phi\in A_1\cap\NR(S)$. Then, $\phi$ cannot be written as a join of atoms of $S$ different from $\phi$. Since a duple positive in $S$ is positive in $M$, by the same proposition, every atom compatible with $M$ is compatible with $S$, so the same holds for $M$. Therefore, $\phi\in\NR(M)$. Since $M$ is complete, $\NWR(M) = \NR(M)$, thus $\phi$ is non-weakly-redundant and, by Corollary \ref{corollary:every_atomization_contains_non-weekly_redundant}, $\phi$ must be in any atomization of $M$. Particularly, $\phi\in A_2$. This shows that $A_1\cap\NR(S)\subseteq A_2\cap\NR(S)$. The proof that $A_1\cap\NR(S)\supseteq A_2\cap\NR(S)$ follows by symmetry. It is now also immediate that $A_1\cap \NWR(S) = A_2\cap\NWR(S)$, since $\NWR(S)\subseteq\NR(S)$.
\end{proof}

We say that a completion of $S$ is the \emph{freest completion} of $S$ if it is freer or as free as every completion of $S$, Definition \ref{definition:freer_or_as_free}; it is unique up to an interpretation-preserving isomorphism, as in the proof of Corollary \ref{corollary:freest_models}. We now show that if we retain all atoms compatible with a model, what we obtain is the freest completion of the model.

\begin{theorem}\label{theorem:freest_completion}
    Let $S$ be a semilattice over a set of constants $C$ with terms $T$. Then, $\Omega(S)$ is an atomization for the freest completion of $S$.
\end{theorem}

\begin{proof}
    Let $M_\Omega$ be the complete semilattice over $C$ spawned by $\Omega(S)$, Proposition \ref{proposition:spawned}, which applies since $\ominus_C\in\Omega(S)$. Since $\Omega(S)$ atomizes $S$, the semilattices $S$ and $M_\Omega$ have the same positive duples over $T$, namely those discriminated by no atom of $\Omega(S)$, so $M_\Omega$ is a completion of $S$ by Proposition \ref{proposition:completion_atomization_base}~(1). Let $M'$ be any completion of $S$ and let $r$ be a duple over $C$, with finite or infinite terms, such that $M'\models r^-$. Some atom $\phi\in\Omega(M')$ discriminates $r$. A duple positive in $S$ is positive in $M'$, by the same proposition, so $\phi$ discriminates no duple positive in $S$, that is, $\phi\in\Omega(S)$, and $M_\Omega\models r^-$. By Proposition \ref{proposition:characterization_of_freedom}, $M_\Omega$ is freer or as free as $M'$. Therefore, $\Omega(S)$ atomizes the freest completion of $S$.
\end{proof}

\begin{corollary}\label{corollary:freest_completion_same_atoms}
    Let $S$ be a semilattice over a set of constants $C$ and let $M$ be the freest completion of $S$. Then, $\Omega(M) = \Omega(S)$.
\end{corollary}

\begin{proof}
    We have seen in Theorem \ref{theorem:freest_completion} that $\Omega(S)$ atomizes $M$. Therefore, $\Omega(S)\subseteq\Omega(M)$. On the other hand, an atom $\phi$ not in $\Omega(S)$ must discriminate some duple $r$ of $S$ such that $S\models r^+$. But since $M$ is a completion of $S$, $M\models r^+$ and, therefore, $\phi\not\in\Omega(M)$. The result now follows.
\end{proof}

We now introduce an illustrative example.

\begin{example}
    Consider $C = \{c_0, c_1, c_2,\dots\}$ and consider $r = (c_0, \odot_{i\ge 1} c_i)$ a duple over $C$. This duple cannot be considered on a finitely generated semilattice, but we can atomize a complete semilattice where that duple is positive by using $A = \{\phi_1,\phi_2,\dots\}\cup\{\phi_{01}, \phi_{02},\dots\}$, where $\phi_i :< \{c_i\}$ and $\phi_{0i} :< \{c_0, c_i\}$. However, in the case of a finitely generated semilattice $S$ over $C$, $A$ is an atomization for the freest model.

    Another atomization for the freest finitely generated semilattice $S$ over $C$ would be $B = \{\phi_0, \phi_1, \phi_2,\dots\}$. In $S$, $\phi_0$ is non-redundant but weakly redundant. That is to say, $\phi_0$ cannot be written as a join of other atoms, but every finitely generated duple discriminated by $\phi_0$ is also discriminated by at least one atom in the set $\{\phi_{01}, \phi_{02},\dots\}$. In fact, the atoms in $B$ are weakly redundant in $S$. Thus, $S$ admits an atomization where all the atoms are non-redundant and weakly-redundant.

    Since $\phi_0$ is an atom compatible with $S$ that is not in the set $A$, nor can it be formed as a union of atoms of $A$, the atomization $A$ is missing some non-redundant atoms of $\Omega(S)$. By Corollary \ref{corollary:freest_completion_same_atoms}, $\phi_0$ is also non-redundant with the freest completion of $S$, which is complete, so $\phi_0$ belongs to every atomization of it. Hence, the complete model generated by $A$ is not the freest completion of $S$. On the other hand, every atom of $\Omega(S)$ is either in $B$ or it can be formed as a union of atoms of $B$. Therefore, by Proposition \ref{proposition:atom_join_right_to_left}, $B$ and $\Omega(S)$ atomize the same model, and in agreement with Theorem \ref{theorem:freest_completion}, $B$ atomizes the freest complete semilattice over $C$.
\end{example}

We now move to results exhibiting how the behaviour of redundant atoms can differ between the finite and infinite cases.

\begin{proposition}\label{proposition:reduntant_not_union_infinite_atom_chain}
    Let $M$ be an atomized semilattice over an infinite set of constants $C$. If $\phi$ is redundant with $M$ but not equal to a union of atoms non-redundant with $M$, then there is an infinite ascending chain $\phi = \phi_0 < \phi_1 < \phi_2 <\dots$ of distinct redundant atoms with $M$.
\end{proposition}

\begin{proof}
    We know that $\phi_0 = \phi$ is redundant, thus it is a union of atoms in $\Omega(M)$. Since $\phi$ is not a union of non-redundant atoms, at least one of these atoms must be redundant and must not be decomposable as a union of non-redundant atoms. Call such atom $\phi_1$. Then $\phi_1$ satisfies the same property, so by applying this process inductively, we eventually obtain the desired chain $\phi = \phi_0 < \phi_1 < \phi_2 <\cdots$.
\end{proof}

This infinite ascending chain of atoms then gives rise to an infinite ascending chain of terms.

\begin{proposition}\label{proposition:ascending_chain_terms}
    Let $M$ be an atomized semilattice over an infinite set of constants $C$ with terms in $T$. If $\phi$ is redundant with $M$ but not equal to a union of atoms non-redundant with $M$, there is an infinite ascending chain of terms of $T$, distinct in $M$, such that $t_1 < t_2 < t_3 < \cdots$.
\end{proposition}

\begin{proof}
    By Proposition \ref{proposition:reduntant_not_union_infinite_atom_chain}, there is an infinite ascending chain of distinct atoms of $M$, $\phi = \phi_0 < \phi_1 < \phi_2 <\cdots$. We then have, for each $i\ge 0$, $U^c(\phi_{i+1})\subsetneq U^c(\phi_i)$.  Since $\phi_0$ may have $C$ as its upper constant segment, we start with $\phi_1$ and, invoking the axiom of choice, select a constant $c_1 \in C \setminus U^c(\phi_{1})$; let $t_1 = c_1$. For atom $\phi_{i + 1}$ select $c_{i + 1} \in  U^c(\phi_{i}) \setminus U^c(\phi_{i + 1})$ and  define $t_{i + 1} = t_{i} \odot c_{i + 1}$.  It follows that $t_{i} \leq t_{i+1}$ in every semilattice and, since $\phi_i <t_{i+1}$ but $\phi_i\nless t_{i}$, then $M \models (t_{i+1} \not\le t_{i})$. Since every term in the sequence $t_i$ is finitely generated $t_i \in T$, and we obtain the desired chain.
\end{proof}

If $M$ is complete, the pinning terms of the atoms in the infinite ascending chain exist and $T_{\phi_1} < T_{\phi_2} < T_{\phi_3} < \cdots$ is an infinite ascending chain of terms of $T$ distinct in $M$.

The following result is now immediate.

\begin{corollary}
    Let $M$ be a semilattice over an infinite set of constants $C$ and terms in $T$. If there is no infinite ascending chain of terms of $T$, distinct in $M$, $t_0 < t_1 < t_2 < \cdots$, every redundant atom with $M$ is a union of atoms non-redundant with $M$.
\end{corollary}

Under the same hypothesis, the picture of finite semilattices is recovered.

\begin{corollary}\label{corollary:ascending_chain_least_atomization}
    Let $M$ be a semilattice over an infinite set of constants $C$ and terms in $T$, and assume there is no infinite ascending chain of terms of $T$, distinct in $M$. Then every atom $\phi\in\Omega(M)$, $\phi\ne\ominus_C$, is pinned by a finite term of $T$. Consequently, $\R(M) = \WR(M)$, $\NR(M) = \NWR(M)$, and $\NR(M)$ is the least atomization of $M$.
\end{corollary}

\begin{proof}
    Let $\phi\in\Omega(M)$, $\phi\ne\ominus_C$, so $C\setminus U^c(\phi)$ is non-empty. Its finite non-empty subsets are terms of $T$. Among their values in $M$ there is a maximal one, $\nu_M(u_0)$; otherwise, we could build an infinite ascending chain of such terms, distinct in $M$. For every $c\in C\setminus U^c(\phi)$, $\nu_M(u_0\cup\{c\})\ge\nu_M(u_0)$, so both are equal by maximality, that is, $M\models c\le u_0$. Hence $u_0$ pins $\phi$.

    By Proposition \ref{proposition:redundancy_and_weak_redundancy}, and by the note after Definition \ref{definition:redundancy} for $\ominus_C$, $\R(M) = \WR(M)$, hence $\NR(M) = \NWR(M)$. By the previous corollary, every atom of $\Omega(M)$ is a union of atoms of $\NR(M)$, so $\NR(M)$ atomizes $M$, Proposition \ref{proposition:atom_join_right_to_left}. By Corollary \ref{corollary:every_atomization_contains_non-weekly_redundant}, it is the least atomization of $M$.
\end{proof}

The hypothesis of Corollary \ref{corollary:ascending_chain_least_atomization} does not depend on $T$, since the values of the terms are the elements of $M$, and, as in its proof, every non-empty set $U$ of constants has a join in $M$, equal to $\nu_M(u_0)$ for some finite $u_0\subseteq U$. Hence $M$ is both finitely generated and complete, and both readings of $M$ have the same compatible atoms and the same atomizations. Indeed, let $(s,t)$ be a duple and let $t_0\subseteq t$ be finite with $\nu_M(t_0)=\nu_M(t)$. An atom compatible with the finitely generated reading is lower than no constant $d$ of $t$, since $(d,t_0)$ is positive, so it discriminates $(s,t)$ if and only if it discriminates $(c,t_0)$ for some $c\in C(s)$, and $(s,t)$ is negative if and only if some such $(c,t_0)$ is. Thus, without an infinite ascending chain of terms, the choice of reading fixed at the beginning of this section makes no difference, unlike in Example \ref{example:infinite_chain}.

We now study infinite ascending chains of atoms in finitely generated semilattices.

\begin{proposition}\label{proposition:limit_atom_finite_generated}
    Let $S$ be a finitely generated atomized semilattice over an infinite set of constants $C$, and let $\Phi$ be a non-empty chain of atoms in $\Omega(S)$ with no greatest element, for instance an infinite ascending chain $\phi_0 < \phi_1 < \phi_2 < \cdots$. If $U = \cap_{\phi\in\Phi} U^c(\phi)\ne \emptyset$, then $\psi:< U$ is in $\Omega(S)$ and is weakly redundant with $S$.
\end{proposition}

\begin{proof}
    Since $U^c(\psi) = U \subseteq U^c(\phi)$ for all $\phi\in\Phi$, and $\Phi$ has no greatest element, $\phi < \psi$ for every $\phi\in\Phi$; in particular $\psi\notin\Phi$. Let $(s,t)$ be a duple in $S$ discriminated by $\psi$. We need to prove that $(s,t)$ is discriminated by some other atom in $\Omega(S)$. We do so by showing that some atom of $\Phi$ must discriminate this duple.

    First, for every $\phi\in\Phi$, $\phi < \psi < s$. On the other hand, we also have $\psi\nless t$. We now show that there must be an atom $\eta\in\Phi$ such that $\eta\nless t$. Assume otherwise that no such $\eta$ exists. Then $C(t)\cap U^c(\phi)$ is non-empty, for every $\phi\in\Phi$. Furthermore, since $S$ is finitely generated, $C(t)$ is finite, and thus the sets $C(t)\cap U^c(\phi)$, $\phi\in\Phi$, are finite and non-empty. They are totally ordered by inclusion, because $\Phi$ is a chain, so there is a least one among them, which is contained in every $U^c(\phi)$ and hence in $U$. Then $C(t)\cap U^c(\psi)$ is also non-empty. But this contradicts that $\psi\nless t$. Finally, $U^c(\psi)\subsetneq U^c(\phi)$ for $\phi\in\Phi$, so $\psi\ne\ominus_C$ and the second condition of weak redundancy holds as well, as shown after Definition \ref{definition:redundancy}. We then have that $\psi$ is weakly redundant with $S$, which by Lemma \ref{lemma:redundant_implies_in_model} implies that $\psi\in\Omega(S)$.
\end{proof}

We can gain some insight on this result by means of Example \ref{example:every_atom_redundant}. In the provided atomization, and for any constant $c\in C$, we can gather infinite ascending chains of atoms so that the intersection of their upper constant segments is $\{c\}$. By Proposition \ref{proposition:limit_atom_finite_generated}, we deduce that, in the finitely generated case, $\phi:<\{c\}$ is an atom weakly redundant with the model. This is another way of showing that in the finitely generated case the provided atomization is an atomization for the freest model, and it further supports that the freest finitely generated model over a set of constants $C$ does not have non-weakly-redundant atoms.

The result also hints that it could be useful to introduce the following concept.

\begin{definition}
    Let $M$ be an atomized semilattice over an infinite set of constants $C$. A limit atom of $M$ is an atom of $\Omega(M)$ whose upper constant segment is equal to the non-empty intersection of the upper constant segments of the atoms of a non-empty chain of atoms with no greatest element, for instance an infinite ascending chain $\phi_0  < \phi_1 < \phi_2 < \cdots$.
\end{definition}

A limit atom may be redundant or non-redundant. However, we can use them to study how we can obtain redundant atoms as unions of atoms.

\begin{theorem}\label{theorem:finitely_generated_join_redundant}
    Let $S$ be a finitely generated atomized semilattice over an infinite set of constants $C$. Every atom redundant with $S$ is equal to a union of atoms non-redundant with $S$.
\end{theorem}

\begin{proof}
    We assume the Axiom of Choice. Let $\kappa$ be the successor cardinal of $2^{|C|}$, an ordinal of cardinality greater than $2^{|C|}$, the cardinal of the set of all possible atoms over $C$, \cite{Kun79}. We will show by contradiction that if there was an atom in $S$ which cannot be written as a union of atoms non-redundant with $S$, then we could find a nested family of different atoms redundant with $S$ indexed on the ordinals $\alpha<\kappa$, thus exceeding $2^{|C|}$. We do so by transfinite recursion, building a family of atoms $\{\phi_\alpha\mid \alpha<\kappa\}$ redundant with $S$ such that if $\beta < \alpha$ then $\phi_\beta < \phi_\alpha$.

    Let $\phi$ be an atom redundant with $S$ and assume that $\phi$ cannot be written as a union of atoms non-redundant with $S$. This implies that there is a constant $c\in U^c(\phi)$ such that if $\psi\in\Omega(S)$ is an atom with $\phi < \psi < c$, then $\psi$ is redundant with $S$. Indeed, first notice that since $\phi$ is redundant, we can necessarily write it as a join of atoms in $\Omega(S)$ different from $\phi$, from which we deduce that for every $c\in U^c(\phi)$ there must be an atom $\psi_c\in\Omega(S)$ with $\phi < \psi_c < c$. If we were able to find one $\psi_c$ that is non-redundant with $S$, for every $c\in U^c(\phi)$, we could write $\phi = \join_c \psi_c$. But then $\phi$ can be written as a join of atoms non-redundant with $S$, contradicting the hypothesis.

    Now let $\phi_0 = \phi$. Let $\alpha$ be an ordinal with $\alpha+1<\kappa$ and assume we have built atoms $\phi_\beta\in\Omega(S)$ redundant with $S$, for all $\beta\le \alpha$, such that $\phi_\beta < \phi_\alpha$ for $\beta < \alpha$, and $\phi_\alpha < c$. Since $\phi_\alpha$ is redundant with $S$, there must be $\psi\in\Omega(S)$ such that $\phi_\alpha < \psi < c$, as for $\phi$ above, and $\psi$ is redundant with $S$ by the choice of $c$, since $\phi\le\phi_\alpha < \psi$. Denote $\phi_{\alpha + 1} = \psi$. We have thus extended the nested sequence of atoms to the successor ordinal of $\alpha$, $\alpha + 1$.

    Assume now that we have a limit ordinal $\lambda<\kappa$ and we have already constructed the atom $\phi_\alpha$ for every $\alpha < \lambda$. Take $U = \cap_{\alpha < \lambda} U^c(\phi_\alpha)$. Then $U\ne\emptyset$ since $c\in U^c(\phi_\alpha)$, for all $\alpha < \lambda$. We can then define $\psi$ to be an atom with $U^c(\psi) = U$, that is, $\psi$ is the limit atom of the atoms we obtained so far. By Proposition \ref{proposition:limit_atom_finite_generated}, applied to the chain $\{\phi_\alpha\mid\alpha<\lambda\}$, which has no greatest element since $\lambda$ is a limit ordinal, $\psi\in\Omega(S)$. Furthermore, $\phi_\alpha < \phi_{\alpha+1}\le\psi$ for all $\alpha < \lambda$, since $\lambda$ is a limit ordinal, and $\psi < c$. In particular, $\phi = \phi_0 < \psi < c$, thus $\psi$ is redundant with $S$ by the choice of $c$. We can then define $\phi_\lambda = \psi$, completing the limit case.

    Consequently, if we assume that an atom $\phi$ cannot be written as a union of atoms non-redundant with $S$, we can build a family of atoms of $S$ indexed on the ordinals $\alpha<\kappa$. These atoms are pairwise different, since $\beta<\alpha$ implies $\phi_\beta<\phi_\alpha$, so the family has cardinality $\kappa$, which exceeds $2^{|C|}$ the cardinality of the set of all atoms over $C$. This is a contradiction. We deduce that $\phi$ must necessarily be equal to a union of atoms non-redundant with $S$.
\end{proof}

\begin{remark}\label{remark:Birkhoff_Frink}
    Theorem \ref{theorem:finitely_generated_join_redundant} can also be proved from the theorem of Birkhoff and Frink that every element of the lattice of subalgebras of an algebra with finitary operations, an algebraic lattice, is a meet of completely meet-irreducible elements, \cite{BirFri48}*{Theorem 7}. The lattice is not $S$ itself but one built from it. Call a set $F\subseteq C$ \emph{closed} if $c\in F$ whenever $S\models (c\le t)$ for a term $t$ with $C(t)\subseteq F$. Since the terms of $S$ are finite, closedness is a closure property of finite character, so the closed sets form under inclusion the lattice of subalgebras of an algebra with finitary operations, \cite{BirFri48}*{Theorem 1}. For an atom $\phi$ write $F_\phi = C\setminus U^c(\phi)$. The atom $\phi$ discriminates $(t,s)$ exactly when $C(s)\subseteq F_\phi$ and $C(t)\not\subseteq F_\phi$, so $\phi\in\Omega(S)$ exactly when $F_\phi$ is closed, and every closed set other than $C$ is $F_\phi$ for a unique $\phi\in\Omega(S)$. Joins of atoms correspond to intersections, $F_{\join_i\phi_i}=\bigcap_i F_{\phi_i}$, so $\phi$ is non-redundant with $S$ exactly when $F_\phi$ is completely meet-irreducible, that is, not an intersection of closed sets different from $F_\phi$. By the theorem of Birkhoff and Frink, $F_\phi$ is an intersection of completely meet-irreducible closed sets $F_{\psi_i}$, so $\phi=\join_i\psi_i$ with every $\psi_i$ non-redundant with $S$. In these terms, Proposition \ref{proposition:limit_atom_finite_generated}, used in the limit step of the proof above, says that the union of a chain of closed sets is closed, \cite{BirFri48}*{Theorem 1}.
\end{remark}

Using this result, we can be a bit more precise on how the atomization of a completion of a finitely generated semilattice $S$ relates to the atomization of $S$, following Theorem \ref{theorem:freest_completion} and Corollary \ref{corollary:freest_completion_same_atoms}.

\begin{proposition}\label{proposition:completion_finitely_generated_NR}
    Every atom of a completion of a finitely generated semilattice $S$ is obtainable as a join of the non-redundant atoms of $S$.
\end{proposition}

\begin{proof}
    Let $M$ be a completion of $S$. Since $M\models\Th_0^+(S)$, Proposition \ref{proposition:completion_atomization_base}~(1), an atom compatible with $M$ discriminates no duple positive in $S$, so $\Omega(M)\subseteq\Omega(S)$. Every atom of $\Omega(S)$ is either non-redundant with $S$ or redundant with $S$, and in the latter case, since $S$ is finitely generated, it is a join of atoms of $\NR(S)$ by Theorem \ref{theorem:finitely_generated_join_redundant}. Hence every atom compatible with $M$ is a join of atoms of $\NR(S)$. In particular, since $\Omega(S)$ atomizes $S$, $S$ itself is atomized by $\NR(S)$, Proposition \ref{proposition:atom_join_right_to_left}.
\end{proof}

We now introduce some characterizations of redundancy.

\begin{proposition}\label{proposition:coinfinite_weakly_redundant}
    Let $M$ be a semilattice over $C$ with terms $T$. Let $\phi$ be a coinfinite atom compatible with $M$. Then either $\phi$ is weakly redundant or there exists $t\in T$ such that for every $c\in C\setminus U^c(\phi)$, $c\le t$.
\end{proposition}

\begin{proof}
    Let $\phi$ be a coinfinite atom compatible with $M$. We will show that if one such term $t\in T$ does not exist, then $\phi$ is weakly redundant. Notice that the pinning term of $\phi$, $T_\phi$, would be one such term. Thus, we can deduce that $T_\phi\not\in T$.

    Thus assume that one such term does not exist, and let $\tilde{M}$ be the freest completion of $M$, Theorem \ref{theorem:freest_completion}. Let $c\in U^c(\phi)$ be any constant. In $\tilde{M}$, we can consider the duple $r = (c, T_\phi)$, which is clearly discriminated by $\phi$. Consequently, $\phi$ is not an atom of $\crossing_r \tilde{M}$. However, $\crossing_r \tilde{M}$ and $M$ must have the same positive theory when restricted to duples in $T$. Indeed, assume that $s=(s_L, s_R)$ is a duple such that $\tilde{M}\models s^-$ and $\crossing_r\tilde{M}\models s^+$. By Proposition \ref{proposition:logical_cosequences}, $\tilde{M}\models T_\phi \le s_R$. From this, we can deduce that $s_R\not\in T$. Indeed, if $s_R\in T$, for every $c\in C(T_\phi) = C\setminus U^c(\phi)$ we would have that $(c,s_R)$ is a duple of $M$. Furthermore, $\tilde{M}\models c \le s_R$, and since $\tilde{M}$ is a completion of $M$, $M\models c \le s_R$. But our hypothesis is that a term in $T$ satisfying this property does not exist.

    Then, since $\crossing_r\tilde{M}$ is complete and has the same positive duples over $T$ as $M$, it is a completion of $M$, Proposition \ref{proposition:completion_atomization_base}~(1), so an atomization of $\crossing_r\tilde{M}$ must also atomize $M$, Proposition \ref{proposition:completion_atomization_base}~(2). Given that $\phi$ discriminates $r$ and thus cannot be an atom of $\crossing_r\tilde{M}$, there must exist an atomization $A$ of $M$ that does not contain $\phi$. Since $\phi$ is compatible with $M$, $A\cup\{\phi\}$ also atomizes $M$, and $A\setminus\{\phi\} = A$. By Theorem \ref{theorem:weakly_redundant_removable}, and since $A\subseteq\Omega(M)$, $\phi$ is weakly redundant.
\end{proof}

Notice that this holds as well if $\phi\ne\ominus_C$ is not coinfinite, but in such case $T_\phi$ always exists, which makes the result trivial. On the other hand, we have the following.

\begin{proposition}
    Let $C$ be a set of constants and let $R$ be a finite set of finite duples over $C$. Let $M = F_C(R)$. If $\phi$ is a coinfinite atom compatible with $M$, it is weakly redundant with $M$.
\end{proposition}

\begin{proof}
By definition, cf.\ Definition \ref{definition:freest_over_duples}, $F_C(R)$ is finitely generated, thus its set of terms $T$ is the set of finite terms over $C$. Let $\phi$ be a coinfinite atom compatible with $M$ and assume that it is not weakly redundant with $M$. Then, there must be a term $t\in T$ such that for every $c\in C\setminus U^c(\phi)$, $c\le t$. Note, however, that $C\setminus U^c(\phi)$ is an infinite set, and $t$ is a finite term. Thus, there is an infinite set of constants not in $C(t)$ which must be below $t$ in $M$. Let $C_R$ be the finite set of constants appearing in the duples of $R$, and choose one such constant $c$ not in $C_R$. Let $N = \{0 < 1\}$ be the two-element semilattice over $C$ with terms $T$ in which $c$ is interpreted as $1$ and every other constant as $0$. Both sides of every duple of $R$ are interpreted as $0$ in $N$, so $N\models R^+$, while $N\models (c\le t)^-$, since $c\notin C(t)$. As $M$ is the freest finitely generated model of $R^+$, it is freer or as free as $N$, and by Proposition \ref{proposition:characterization_of_freedom}, $M\models (c\le t)^-$, a contradiction. Thus $\phi$ must be weakly redundant with $M$.
\end{proof}

However, it is still the case that finitely generated models can have coinfinite non-weakly-redundant atoms.

\begin{example}
    Let $C = \{c_i\mid i\in\mathbb{N}\}$. Consider $M$ to be the infinite descending chain $c_0 > c_1 > c_2 >\cdots$. As a semilattice, it is finitely generated over $C$, thus we consider its set of terms $T$ to be the set of finite terms over $C$. It is also easy to check that $M$ is atomized by $A = \{\phi_i :<\{c_0, c_1, \dots, c_i\}\colon i\in\mathbb{N}\}$.

    These are, in fact, the only atoms compatible with the model other than $\ominus_C$. Indeed, assume that $\phi$ is an atom compatible with $M$ and that $c_i\in U^c(\phi)$. If there is $j < i$ such that $c_j\not\in U^c(\phi)$, $\phi$ discriminates $(c_i, c_j)$, but $M\models c_j > c_i$. Thus, $c_j\in U^c(\phi)$, for all $j\le i$. Consequently, $U^c(\phi)$ is either $\{c_0,\dots,c_i\}$ for some $i$, or $C$, and $\Omega(M) = A\cup\{\ominus_C\}$. Since $M$ has no least element, $\ominus_C = \join_i\phi_i$ is redundant. We can now say that each of the $\phi_i$ is non-weakly-redundant with $M$, as $\phi_i$ is the only atom in $\Omega(M)$ discriminating $(c_i, c_{i+1})$. Nonetheless, $\phi_i$ is also coinfinite.

    Notice, however, that given $\phi_i$, $C\setminus U^c(\phi_i) = \{c_j\mid j > i\}$. Then, $c_{i+1}$ is a finite term such that $c\le c_{i+1}$, for every $c\in C\setminus U^c(\phi_i)$. Thus, $\phi_i$ being non-weakly-redundant with $M$ is consistent with Proposition \ref{proposition:coinfinite_weakly_redundant}.
\end{example}

We can also prove the following result.

\begin{proposition}\label{proposition:atom_classification_finitely_generated}
    Let $M$ be any semilattice over $C$.
    \begin{enumerate}
        \item If a cofinite atom in $\Omega(M)$ is weakly redundant with $M$, then it is redundant.
        \item An atom weakly redundant with $M$ that is non-redundant with $M$ is coinfinite.
    \end{enumerate}
\end{proposition}

\begin{proof}
    To prove (1), assume that $\phi\in\Omega(M)$ is cofinite. If $\phi=\ominus_C$, it is redundant if and only if it is weakly redundant, as shown after Definition \ref{definition:redundancy}. Otherwise, its pinning term $T_\phi$ exists, thus the result follows from Proposition \ref{proposition:redundancy_and_weak_redundancy} (2). The second point now follows immediately.
\end{proof}

\begin{corollary}
    Let $M$ be a semilattice over a set of constants $C$. If $\phi\in\Omega(M)$ is non-redundant and cofinite, then $\phi\in A$ for any atomization $A$ of $M$.
\end{corollary}

\begin{proof}
    As a consequence of Proposition \ref{proposition:atom_classification_finitely_generated} (1), a non-redundant cofinite atom is non-weakly-redundant. The result is then a consequence of Corollary \ref{corollary:every_atomization_contains_non-weekly_redundant}.
\end{proof}

\section*{Acknowledgements}

This work was supported by Champalimaud Research, Lisbon, and by the H2020 ICT48 project ALMA: Human Centric Algebraic Machine Learning (grant 952091).

\end{document}